\documentclass[11pt,a4paper]{article}
\usepackage[latin1]{inputenc}
\usepackage[T1]{fontenc}
\usepackage{amsmath}
\usepackage{amsfonts}
\usepackage{amssymb}
\usepackage{theorem}
\usepackage{enumerate}
\usepackage{mathrsfs}
\usepackage{dsfont}
\usepackage{stmaryrd}
\usepackage{cancel}

\renewcommand{\geq}{\geqslant}
\renewcommand{\leq}{\leqslant}

\newcommand{\reff}[1]{(\ref{#1})}

\newtheorem{thm}{Theorem}
\newtheorem{lem}{Lemma}
\newtheorem{cor}[lem]{Corollary}
\newtheorem{prop}[lem]{Proposition}
\newtheorem{rem}{Remark}[section]
\newtheorem{Example}{Example}[section]

 \newtheorem{condition}[thm]{Condition}

{\theorembodyfont{\rmfamily}\newtheorem{comment}{Comment}[section]}

 \def\blemma{\begin{lem}}\def\elemma{\end{lem}}
 \def\bproposition{\begin{prop}}\def\eproposition{\end{prop}}
 \def\btheorem{\begin{thm}}\def\etheorem{\end{thm}}
 \def\bcorollary{\begin{cor}}\def\ecorollary{\end{cor}}
 \def\bremark{\begin{rem}}\def\eremark{\end{rem}}
 \def\bcondition{\begin{condition}}\def\econdition{\end{condition}}

 \def\benumerate{\begin{enumerate}}\def\eenumerate{\end{enumerate}}
 \def\bitemize{\begin{itemize}}\def\eitemize{\end{itemize}}

 \def\beqlb{\begin{eqnarray}}\def\eeqlb{\end{eqnarray}}
 \def\beqnn{\begin{eqnarray*}}\def\eeqnn{\end{eqnarray*}}

\newcommand{\E}{{\bf E}}
\renewcommand{\P}{{\bf P}}

\newcommand{\ddr}{\mathrm{d}}
\newcommand{\R}{\mathbb{R}}
\newcommand{\bbN}{\mathbb{N}}

\newcommand{\bbQ}{\mathbb{Q}}

\newcommand{\un}{\mathbf{1}}

\newcommand{\cC}{\mathcal{C}}

\newcommand{\cF}{\mathcal{F}}
\newcommand{\cG}{\mathcal{G}}
\newcommand{\cI}{\mathcal{I}}
\newcommand{\cJ}{\mathcal{J}}
\newcommand{\cL}{\mathcal{L}}
\newcommand{\cM}{\mathcal{M}}
\newcommand{\cN}{\mathcal{N}}

\newcommand{\cT}{\mathcal{T}}

\newcommand{\N}{\bold{N}}

\newcommand{\veps}{\varepsilon}
\newcommand{\noi}{\noindent}

\newcommand{\lb}{\ell^b}
\newcommand{\la}{\ell^a}

\newcommand{\D}{\delta}

\newcommand{\talpha}{\tilde{\alpha}}

\newcommand{\lgeo}{[{\!  } [}
\newcommand{\rgeo}{] {\!  } ]}
\newcommand{\cH}{\mathscr{H}}
\newcommand{\ccB}{\mathscr{B}}
\def\cqfd{\hfill $\blacksquare$ }

\newcommand{\br}{\textbf{r}}
\newcommand{\beps}{{\boldsymbol \varepsilon}}
\newcommand{\Rp}{\br^{(p)}}
\newcommand{\bepsp}{\beps^{(p)}}

\newcommand{\bC}{\mathbf{C}}
\newcommand{\bV}{\mathbf{V}}

\newcommand{\supp}{\mathtt{Supp}}

\title{Uniform Hausdorff measure of the level sets \\of the Brownian tree.}

\author{ Xan \textsc{Duhalde}\thanks{\textbf{Institution}: PRES Sorbonne Universit\'es, UPMC Universit\'e Paris 06,  LPMA (UMR 7599). \textbf{Postal address}: LPMA, Bo\^ite courrier 188, 4 place Jussieu, 75252 Paris Cedex 05, FRANCE. \textbf{Email}: xan.duhalde@upmc.fr}}
\vspace{2mm}

\begin{document}
\maketitle
\begin{abstract}
Let $(\cT,d)$ be the random real tree with root $\rho$ coded by a 
Brownian excursion. So $(\cT,d)$ is (up to normalisation) Aldous CRT \cite{AldousI} (see Le Gall \cite{LG91}). The $a$-level set of $\cT$ is the set $\cT(a)$ of all points in $\cT$ that are at distance $a$ from the root. We know from Duquesne and Le Gall \cite{DuLG06} that for any fixed $a\in (0, \infty)$,  the measure $\ell^a$ that is induced on $\cT(a)$ by the local time at $a$ of the Brownian excursion, is equal, up to a multiplicative constant, to the Hausdorff measure in $\cT$ with gauge function $g(r)= r \log\log1/r$,  restricted to $\cT(a)$.   
As suggested by a result due to Perkins \cite{Per88,Per89} for super-Brownian motion, we prove in this paper  
a more precise statement that holds almost surely 
uniformly in $a$, and we specify the multiplicative constant. Namely,  
we prove that almost surely for any $a\in (0, \infty)$, $\ell^a(\cdot) = \frac{1}{2} \cH_g (\, \cdot \, \cap \cT(a))$, where $\cH_g$ stands for the $g$-Hausdorff measure.

\smallskip


\smallskip


\end{abstract}

\section{Introduction.}

The Continuum Random Tree  was introduced by Aldous \cite{AldousI} as a random compact  metric space $(\cT_{\small{1}},d,\textbf{m}_{\small{1}}) $, endowed with a mass measure $\textbf{m}_{\small{1}}$ such that  almost surely $\textbf{m}_{\small{1}}(\cT_{\small{1}})=1$.   It appears as the scaling limit of a large class of discrete models of  random trees,  and can be alternatively encoded by a normalised Brownian excursion (see Le Gall \cite{LG91}). This encoding procedure will be the viewpoint of  the present  paper, but  for the sake of simplicity, we will \textit{not} ask the total mass to be equal to one. Instead,  we  work on the tree encoded by a Brownian excursion $(e_t, t\geq 0)$, under  its excursion measure $\N$. Let us mention that our result remains true for the CRT. 

\par The Brownian tree has a distinguished vertex $\rho$  called the root, so it makes sense to define, for all $a\in(0,\infty)$ the $a$-level set 
$\cT(a)=\left\{\sigma\in\cT : d(\rho,\sigma)=a\right\}$. Moreover, one can define the collection of measures $\left(\la(\ddr \sigma), \sigma\in\cT, a\in(0,\infty)\right)$,
as the image of the local times on the levels of the excursion. Those measures are called local time measures. Indeed, $\N$-a.e. for all $a\in(0,\infty)$, the topological support of $\la$ is included in $\cT(a)$. Duquesne and Le Gall \cite{DuLG06} showed that for a fixed level $a$, one has  
\begin{equation}\label{resultatfixe}
\textrm{$\N$-a.e.}\quad\ell^a(\cdot) = c\cH_g (\, \cdot \, \cap \cT(a)),
\end{equation}
where $\cH_g$ stands for the Hausdorff measure associated with the gauge function $g(r)=r \log\log 1/r$ and $c\in(0,\infty)$ is a multiplicative constant. In this paper, we prove that $c=\frac{1}{2}$ and that the result holds $\N$-a.e.~simultaneously for all levels $a$. Let us mention that the value $\frac{1}{2}$ depends on the normalisation chosen for the excursion measure $\N$. A similar  result has been obtained by Perkins \cite{Per88,Per89} for Super Brownian Motion.  Briefly, let $(Z_a, a\geq 0)$ a version of this measure-valued process on $\R^d$, defined on $(\Omega,\cF,\P)$.  Perkins proves that if the dimension $d$ of the space is such that  $d\geq 3$ (which corresponds to the supercritical dimension case),  there exists two constants $c_d, C_d$ in $(0,\infty)$, only depending on $d$ such that  the following holds 
\begin{equation}\label{resultatPerkins}
\textrm{$\P$-a.s.}\quad \forall a\in(0,\infty) \quad
c_d\cH_g\left(\cdot\cap \supp(Z_a)\right)\leq Z_a\left(\cdot\right)\leq C_d\cH_g\left(\cdot\cap \supp(Z_a)\right),
\end{equation} where $\supp(Z_a)$ is the topological support of the measure $Z_a$ and $\cH_g$ is the Hausdorff measure associated to the gauge function $g(r)=r^2\log\log 1/r$. In this paper, we use the ideas and techniques of \cite{Per88,Per89} to get a result similar to (\ref{resultatPerkins}), an equality being accessible in the setting of trees.\\

\par Before stating formally our result, let us recall precisely basic facts.
A metric  space $(T,d)$ is a real tree if and only if the following two properties hold for any $\sigma_1, \sigma_2$ in $T$ : 

\begin{enumerate}
\item[(i)] There is a unique
isometric map
$f_{\sigma_1,\sigma_2}$ from $[0,d(\sigma_1,\sigma_2)]$ into $T$ such
that $f_{\sigma_1,\sigma_2}(0)=\sigma_1$ and $f_{\sigma_1,\sigma_2}(
d(\sigma_1,\sigma_2))=\sigma_2$. We set $\llbracket \sigma_1,\sigma_2\rrbracket=f_{\sigma_1,\sigma_2}\left([0,d(\sigma_1,\sigma_2)]\right)$ that is the geodesic path  joining $\sigma_1$ and $\sigma_2$.
\item[(ii)] If $q$ is a continuous injective map from $[0,1]$ into
$T$, such that $q(0)=\sigma_1$ and $q(1)=\sigma_2$, we have
$$q([0,1])=f_{\sigma_1,\sigma_2}([0,d(\sigma_1,\sigma_2)]).$$
\end{enumerate}
If $\sigma_1\in\llbracket\rho,\sigma_2\rrbracket$, we will say that $\sigma_1$ is an \textit{ancestor} of $\sigma_2$ ($\sigma_2$ is a \textit{descendant} of $\sigma_1$).\\

Real trees can be derived from continuous functions that represent their contour functions. Namely, let us consider a (deterministic) excursion  $e$, that is to say a continuous function for which  there exists $\zeta\in(0,\infty)$ such that : $\forall t\geq \zeta, e(0)=e(t)=0$, and  $\forall t\in(0,\zeta), e(t)>0$. A real tree $T$ can be associated with $e$ in the following way. For $s,t\in[0,\zeta]$, we set
$$ d(s,t)=e(s)+e(t)-2\inf_{r\in[s\wedge t,s\vee t]}e(r).$$
It is easy to see that $d$ is a pseudo-distance on $[0,\zeta]$. Defining the equivalence relation 
$s\sim t$ iff $d(s,t)=0$, one can set 
\begin{equation}\label{quotient}
T=[0,\zeta]/ \sim.
\end{equation}
The function $d$ induces a distance on the quotient set $T$. For a fixed excursion $e$, let 
\begin{equation}\label{defpe}
p: [0,\zeta] \longrightarrow (T,d)
\end{equation}
be the canonical projection. Clearly $p$ is continuous, which implies that $(T,d)$ is a compact metric space. Moreover, it can be shown (see \cite{DuLG05} for a proof) that $(T,d)$ is a real tree\\
We take  $\rho=p(0)$ as the root of $T$. For all $a\in(0,\infty)$, the $a$-level set $T(a)=\left\{\sigma \in T : d(\rho,\sigma)=a\right\}$
is the image by $p$ of the set $\{ t\in[0,\zeta] : e(t)=a\}$.
The total \textit{height} of the tree is defined by
\begin{equation}\label{defheight}
h(T)=\sup\left\{d(\rho,\sigma); \sigma\in T\right\}.
\end{equation}

We define the  \textit{Brownian tree} as  the metric space $(\cT,d)$ coded by the Brownian excursion. More precisely,  let $(\Omega,\cF,\P)$ a probability space, large enough to carry all the random variables we need. We consider on that space a process $(X_t, t\in [0,\infty))$ such that $(\frac{1}{\sqrt{2}}X_t, t\in [0,\infty))$ is a standard real-valued Brownian motion (the choice of the normalizing constant $\sqrt{2}$ is explained below). Let us set $\underline{X}_t =\inf_{s\in [0, t]} X_s$. Then, the reflected process $X-\underline{X}$ is a strong Markov process, and the state $0$ is instantaneous in $(0,\infty)$ and recurrent (see \cite{Bebook}, chapter VI). We denote by $\N$ the excursion measure associated with the local time $-\underline{X}$; $\N$  is a sigma-finite measure on the space of continuous functions on $[0,\infty)$, denoted $\bC^0$ in this work. More precisely, let $\bigcup_{j\in \cJ} (l_j, r_j)= \{ t>0: X_t-\underline{X}_t >0\}$ be the excursion intervals of the reflected process, and for all $j\in \cJ$, we set $e_j (s)= X_{(l_j +s)\wedge d_j} -  \underline{X}_{l_j}$, $s\in [0, \infty)$. Then,
$$ \mathcal{M}(\ddr t,\ddr e)= \sum_{j\in \cJ} \delta_{(-\underline{X}_{l_j}, e_j)}$$
is a Poisson point measure on $[0,\infty)\times\bC^0$ of intensity $\ddr t \N(\ddr e)$. Let us recall that  the two processes $\left(\lvert X_t\rvert,2L_t\right)_{t\geq 0}$ and $\left(X_t-\underline{X}_t,-\underline{X}_t\right)_{t\ge 0}$ have the same law under $\P$ by a celebrated result of L\'evy (see Blumenthal \cite{Blubook}, Th. II 2.2) where  the process $(L_t, t\geq 0)$ is defined by the approximation $L_t=\lim\limits_{\veps\to 0}(2\veps)^{\!-1}\int_0^t\un_{\{\lvert X_s\rvert\leq \veps\}}\ddr s$ that holds uniformly in $t$ on compact subsets of $[0,\infty)$. \\
 
We shall denote by $(e_t, t\geq 0)$ the canonical process on $\bC^0$. Under $\N$, it is a strong Markov process, with transition kernel of the original process $X$ killed when it hits $0$ (see \cite{Blubook} III 3(f)). The following properties hold for the process $\N$-a.e. : there exists  a unique real $\zeta\in(0,\infty)$ such that $\forall t\in (0,\zeta), e(t)>0$, and $\forall t\in[\zeta,\infty), e(t)=e(0)=0$. Moreover, with our normalization, one has (see \cite{Blubook} IV 1.1)
\begin{equation}\label{loizeta}
\forall\lambda\in[0,\infty), \mathbf{N} (1-e^{-\lambda \zeta})= \sqrt{\lambda}\quad {\rm and}\quad
\mathbf{N} (\zeta \in \ddr r)= \frac{r^{-3/2} }{2\sqrt{\pi}}\ddr r.
\end{equation}
One can  show that $\N\left(\cdot\mid\zeta\in[1-\veps,1+\veps]\right)$ converges when $\veps$ goes to $0$, towards a probability measure that is denoted by $\N(\cdot \mid \zeta=1)$. It can be seen as the law of the excursion of $X-\underline{X}$ conditioned to have length one. The tree encoded by $e$ under $\N(\cdot \mid \zeta=1)$ is the CRT defined in \cite{AldousI}. 
The choice of the normalising constant $\sqrt{2}$ is explained by the following. 
Let $\tau_n$ be uniformly distributed as the set of rooted planar trees with $n$ vertices. We view it as a real tree, the edges of $\tau_n$ being intervals of length one, and we denote by $(\tau_n,d_n)$ the resulting metric space. Denote by $(C_{t}^{(n)}, t\in[0,2(n\!-\!1)])$ its contour function that is (informally) defined as follows. We let a particle explore the planar tree at speed one, from the left to the right, beginning  at the root. We set $C_t^{(n)}$ as the distance from the root of the particle at time $t$. It can be shown (see \cite{LGCornell} Th. 1.17) that $(C_{t}^{(n)}, t\in[0,2(n\!-\!1)])$ has the law of a simple random walk conditioned to be positive on $[1,2(n\!-\!1)-1]$ and null at $2(n\!-\!1)$. Using Donsker invariance principle, the rescaled contour function $(n^{-1/2}C_{2(n\!-\!1)t}^{(n)}, t\in[0,1])$ converges in law towards the law of $(e_t, t\in[0,1])$ under $\N(\cdot \mid \zeta=1)$. In terms of  trees,  $(\tau_n, n^{-\!1/2}d_n)$ converges towards the CRT, that is the tree  $(\cT_{\small{1}},d) $ coded by $e$ under $\N(\cdot \mid \zeta=1)$. The latter convergence can be stated using the distance of Gromov-Hausdorff (see Evans, Pitman, Winter \cite{EPW06}).

Recalling definition (\ref{defheight}), we get from \cite{Blubook} IV 1.1 that with our normalization,

\begin{equation}\label{rephT}
\forall a\in(0,\infty)\quad \N\Big(\sup\limits_{t\in[0,\zeta]} e_t>a\Big)=\N\Big(h(\cT)>a\Big)=\frac{1}{a}.
\end{equation}
In the paper, for $a\in(0,\infty)$ we shall use the probability measure, 
\begin{equation}\label{notationNa}
\N_a=\N\left(\cdot\mid h(\cT)>a\right)=a\N\left(\cdot\un_{\{h(\cT)>a\}}\right).
\end{equation}
Recall that the $a$-level set of the Brownian tree is  defined by 
\begin{equation}
\cT(a)=\left\{\sigma \in \cT : d(\rho,\sigma)=a\right\}.
\end{equation}
As a consequence of Trotter's theorem on the regularity of Brownian local time (\cite{Blubook} sec VI.3)  there exists a $[0, \infty)$-valued process $(L^a_t)_{a, t\in [0, \infty)}$ such that $\N$-a.e.~the following holds true: 
\begin{itemize}
\item $(a, t) \mapsto L^a_t $ is continuous, 
\item for all $a\in [0, \infty)$, $t\mapsto L^a_t$ is non-decreasing, 
\item for all $a\in [0, \infty)$, for all $t\in [0, \infty)$ and for all $b\in (0, \infty)$, 
\begin{equation}
\label{loctimeexc}
\lim_{\veps \to 0} \N \Big( \un_{\{ \sup e >b \}} \sup_{0\leq s \leq t\wedge \zeta} \Big| \frac{1}{\veps} \int_0^s\!\!\! \un_{\{ a-\veps < e(u) \leq a \}} \ddr u -L^a_s \Big|  \Big)=0 \; .
\end{equation}
\end{itemize}
We refer to \cite{DuLG}, Proposition 1.3.3. for details in a more general setting.\\

The image by the projection $p:[0,\zeta]\to\cT$ of those local times defines the collection of local time measures on the tree, $\left(\la(\ddr \sigma), \sigma\in\cT, a\in(0,\infty)\right)$. More precisely,
\begin{equation}\label{defella}
\textrm{$\N$-a.e.} \quad {\rm for\ all\ }f: \cT\overset{{\rm meas.}}{\to} [0,\infty)\quad \forall a\in(0,\infty) \quad \int_\cT f(\sigma ) \ell^a (\ddr \sigma)= \int_0^\zeta f(p(t)) \ddr L^a_t.
\end{equation}
See \cite{DuLG05}, Th. 4.2 for an intrinsic definition  of the measure $\la$ (for fixed $a$).  
Let $\cG_a$ the $\sigma$-field generated by the excursion below level $a$ (formal definitions and details on what follows are given in section \ref{basicfactsBrownian}). The approximation (\ref{loctimeexc}) entails that for fixed $a$,  $\la(\cT)=L^a_\zeta$ is $\cG_a$ measurable. Moreover,  the Ray-Knight theorem (\cite{Blubook} VI 2.10) entails that under $\N_a(\cdot)$ conditionally  on $\cG_a$, the process $\left(\ell^{a+a^\prime}(\cT), a^\prime\geq 0\right)$ is a Feller diffusion started at $\la(\cT)$. In particular, one has 
\begin{equation}\label{Laplaceella}
\forall a,\lambda\in(0,\infty)\quad\N\left[1-e^{-\lambda \ell^a(\cT)}\right]=\frac{\lambda}{1+a\lambda},
\end{equation} 
which implies that under $\N_a$, $\la(\cT)$ is exponentially distributed with mean $a$. The regularity of $a\mapsto\la(\cT)$ is extended by Duquesne and Le Gall \cite{DuLG05} : they prove that $\N$-a.e. the process $a\mapsto \la$ is continuous for the weak topology of measures. In the same work, the  topological support of the level set measures is precised as follows. A vertex $\sigma\in\cT$ is called an extinction point if there exists $\veps\in(0,\infty)$ such that $d(\rho,\sigma)=\sup\{d(\rho,\tau), \tau\in B(\sigma,\veps)\}$, where $B(\sigma,\veps)$ is the open ball in $\cT$ with center $\sigma$ and radius $\veps$. For $s\in[0,\zeta]$, the vertex $p(s)\in\cT$ is an extinction time iff $s\in[0,\zeta]$ is a local maximum of $e$. As a consequence, the set of all extinction points, denoted $\mathscr{E}$, is countable. Let us denote $\supp(\mu)$ for the topological support of the measure $\mu$. The result states that 
\begin{equation}\label{supportla}
\textrm{$\N$-a.e.}\quad
\forall a\in(0,\infty)\setminus\mathscr{E}, \supp(\la)=\cT(a), \quad \textrm{and}\quad
\forall a\in\mathscr{Ext}, \supp(\la)=\cT(a)\setminus \{\sigma_a\},
\end{equation}
where $\sigma_a$ is the (unique) extinction point at level $a$ (see Perkins \cite{Per90} for previous results on Super-Brownian motion). \\

\par  Let us briefly introduce the construction of  the Hausdorff measure. We set the \textit{gauge} function  $g$ as
\begin{equation}
g(r)=r\log\log 1/r,\quad r\in(0,e^{\!-\!1}).
\end{equation}
In all the paper it will be assumed implicitly that $g(r)$ is considered only for $r\in(0,e^{\!-\!1})$. On that interval, $g$ is an increasing continuous function. For any subset $A$ of $\cT$, one can define
\begin{equation}\label{defHausdorff}
\cH_g(A)=\lim\limits_{\veps\to 0}\ \inf\left\{\sum\limits_{i\in\bbN} g\left({\rm diam}(E_i)\right) \ ; A\subset \bigcup\limits_{i\in\bbN} E_i, {\rm diam}(E_i)<\veps\right\}.
\end{equation} 
Standard results on Hausdorff measures (see e.g.  \cite{Rogersbook}) 
ensure that $\cH_g$ defines a Borel-regular  outer measure on $\cT$ called the $g$-Hausdorff measure on $\cT$.  The main result of the paper is the following.

\begin{thm}\label{mainth}
Let $\cT$ be the Brownian tree, that is the tree encoded by the  excursion $e$ under $\N$. Let $(\la(\ddr\sigma), \sigma\in\cT, a\in(0,\infty))$ the collection of local time measures and $\cH_g$ the $g$-Hausdorff measure on $\cT$, where $g(r)=r\log\log 1/r$. Then, the following holds :  
\begin{equation}\label{mainresult}
\textrm{$\N$-a.e.}\quad\forall a\in (0,\infty)\quad\la(\cdot)={1\over 2}\cH_g\left(\cdot\cap\cT(a)\right).
\end{equation}
\end{thm}

\begin{comment}
Thanks to the scaling properties of the Brownian excursion, one can derive from Theorem~\ref{mainth} a similar statement for the tree coded by $e$ under $\N(\cdot \mid \zeta=1)$, that is Aldous CRT. 
\end{comment}

\begin{comment}
Our result seems close to a theorem of Perkins \cite{Per81} on linear Brownian motion.  Let $(L^a_t, t\geq 0, a\in\R)$ be the bi-continuous version of the local times for the process $(X_t, t\geq 0)$ defined above. Those local times are given by an approximation of the type of  (\ref{loctimeexc}). Perkins proves that almost surely, uniformly in $a$, one has $L_t^a=\cH_g(\{s\in[0,t] : X_s=x\})$, where $\cH_g$ stands for the Hausdorff measure on the line associated with the gauge  $g(r)=\sqrt{r\log\log 1/r}$ (the result for fixed $a$ had been obtain by Taylor and Wendel in \cite{TW66}). The Brownian tree being coded by the Brownian excursion, everything happens as if the projection mapping $p : [0,\zeta]\to \cT$ is $1/2$-Hölder and induces a strong "doubling", such that the entire gauge function is squared. Nevertheless, we don't see how to derive  our result from \cite{Per81}. 
\end{comment}

The paper is organised as follows. In section \ref{secgeom}, we state some deterministic facts on the geometry of the level sets for a real tree. In particular, we provide two comparison lemmas with respect to Hausdorff measure on real trees. The second one, that is specific to our setting, seems new to us. In section \ref{secBtree}, we recall basic facts on the Brownian tree and we establish some technical estimates. Section \ref{preuvetheoreme} is devoted to the proof of Theorem \ref{mainth}. As a first step, we prove Theorem \ref{upper}, which gives an upper bound for the local time measures. To that end, we need to control the total mass of the balls that are "too large". The second step is the proof of Theorem \ref{BP}, which requires a control of the number of balls that are "too small". Let us mention again  that our strategy and many ideas in this work were borrowed from \cite{Per88,Per89}.

\bigskip

\noi \textbf{Acknowledgments}. I would like to thank my advisor Thomas Duquesne for introducing this problem, as well  as for his help and the many improvements he suggested.

\section{Geometric properties of the level sets of real trees.}\label{secgeom}

\subsection{The balls of the level sets of real trees.}\label{subsecballs}

Let $(T,d,\rho)$ be a compact rooted real tree as defined in the introduction. Recall that for any $\sigma, \sigma^\prime \in T$, $\lgeo \sigma , \sigma^\prime \rgeo$ stands for the unique 
geodesic path joining $\sigma$ to $\sigma^\prime$. We shall view $T$ as a family tree whose ancestor is the root $\rho$ and we then denote by $\sigma \wedge \sigma^\prime$ the most recent common ancestor of $\sigma$ and $\sigma^\prime$ that is formally defined by 
$$ \lgeo \rho , \sigma \wedge \sigma^\prime \rgeo = \lgeo \rho , \sigma  \rgeo \cap \lgeo \rho ,  \sigma^\prime \rgeo \; .$$
Observe that 
\begin{equation}
\label{hauteurdist}
\forall \sigma, \sigma^\prime \in T , \quad d(\sigma, \sigma^\prime)= d(\rho, \sigma) + d(\rho, \sigma^\prime)-2d(\rho, \sigma \wedge \sigma^\prime) \; .
\end{equation}
Let $a\in [0, \infty)$. Recall that the $a$-level set of $T$ is given by 
$$T(a)= \big\{ \sigma \in T: d(\rho, \sigma )= a \big\}. $$

\paragraph{Subtrees above level $b$.} Let $b\in [0, \infty)$ and 
denote by $(T^{o,b}_j)_{j\in \cJ_b}$ the connected components of the open set 
$\{ \sigma \in T: d(\rho, \sigma) >b\}$:  
$$ \bigcup_{j\in \cJ_b} T^{o,b}_j = \big\{ \sigma \in T: d(\rho, \sigma) >b \big\} \; .$$
Then for any $j\in \cJ_b$, there exists a unique point $\sigma_j\in T(b)$ such that $T^b_j:= T^{o,b}_j \cup \{ \sigma_j\}$ is the closure of $T^{o,b}_j$ in $T$. 
Note that $(T_j^b, d, \sigma_j)$ is a compact rooted real  tree and that 
$$ \forall j\in \cJ_b, \; \forall \sigma \in T^b_j, \quad \sigma_j \in \lgeo \rho, \sigma \rgeo   \; .$$

\paragraph{Open balls in $T(a)$.} Recall that $B(\sigma, r)$ stands for the open ball in $T$ with center $\sigma$ and radius $r$. We shall also denote by $\Gamma (\sigma, r)$ the open ball with center $\sigma$ and radius $r$ \textit{in the level set of $\sigma$}, namely 
\begin{equation}
\label{Taboule}
\Gamma (\sigma, r)= B(\sigma, r) \cap T(a) \, , \quad \textrm{where $a=d(\rho, \sigma)$.}
\end{equation}
If $\sigma \in T(a)$, then we call $\Gamma (\sigma, r)$ a $T(a)$-ball with radius $r$; we denote by $\ccB_{a,r}$ the set of all the $T(a)$-balls with radius $r$: 
\begin{equation}
\label{Tabouleset}
\ccB_{a,r}= \big\{ \Gamma (\sigma, r) ; \sigma \in T(a)  \big\} \; .
\end{equation}
The following proposition provides the geometric properties of $T(a)$-balls that we shall use. 

\begin{prop}
\label{propGamm} Let $(T,d, \rho)$ be a compact rooted real tree. Let $a, r\in (0, \infty)$ be such that $a\geq r/2$. Then, the number of $T(a)$-balls with radius $r$ is finite. We set 
\begin{equation}
\label{defZar}
 Z_{a,r} = \# \ccB_{a, r} \quad \textrm{and} \quad  \big\{\,  \Gamma_i, \; 1 \! \leq\!  i \! \leq\!  Z_{a,r}\big\}= \ccB_{a,r}. 
\end{equation} 
Then , the following holds true. 
\begin{itemize}
\item[(i)] Set $b\! =\! a\! -\! \frac{_1}{^2}r$. Then, there are $Z_{a,r}$ distinct subtrees above $b$ 
denoted by $(T^{b}_{j_i}, d, \sigma_{j_i})$, $j_i \in \cJ_b$, $1\leq i\leq Z_{a,r}$ such that  
$$ \Gamma_i = T(a) \cap T^{b}_{j_i} = \{ \sigma^\prime \in T^{b}_{j_i}: d(\sigma_{j_i}, \sigma^\prime)=r/2 \big\} \; .$$ 
Thus, the $T(a)$-balls with radius $r$ are pairwise disjoint. 
\item[(ii)] For all $\sigma \!  \in \! T(a)$, one has ${\rm diam} (\Gamma (\sigma, r)) \! \leq \! r$. If furthermore
$r \! \in \! (0, 2a)$, then 
${\rm diam} (\Gamma (\sigma, r)) \! < \! r$ and 
\begin{equation}
\label{diamboule}
 \forall r^\prime \in  \big( {\rm diam} (\Gamma(\sigma, r)), r \big) \quad \Gamma (\sigma, r^\prime)= \Gamma (\sigma, r)\; .
\end{equation} 
Therefore, the set of all $T(a)$-balls is countable.  
\item[(iii)] Two $T(a)$-balls are either contained one in the other or disjoint. Namely, for all 
$r^\prime<r$ and all $\sigma, \sigma^\prime \in T(a)$, either 
$\Gamma(\sigma^\prime, r^\prime) \subset   \Gamma(\sigma, r)$ or $\Gamma(\sigma^\prime, r^\prime) \cap  \Gamma(\sigma, r)=\emptyset$. 
\end{itemize}
\end{prop}

\noindent{\bf Proof.} Let us prove $(i)$. Let $\sigma , \sigma^\prime \in T(a)$ and set $b\! =\! a\! -\! \frac{_1}{^2}r$. By (\ref{hauteurdist}), $d(\sigma, \sigma^\prime)= 2a-2d(\rho, \sigma\wedge \sigma^\prime)$. Thus, $d(\sigma, \sigma^\prime) <r$ iff $d(\rho, \sigma\wedge \sigma^\prime) >b$. 
Let $j\in \cJ_b$ be such that $\sigma \in T^b_j$; namely, $T^b_j$ is the unique subtree above $b$ containing $\sigma$ and 
$\sigma_j$ is the unique point $\gamma \in \lgeo \rho , \sigma \rgeo$ such 
that $d(\rho, \gamma)= b$. Now observe that for all $\sigma^\prime \in T(a)$, 
$$ d( \rho , \sigma \wedge \sigma^\prime) > b \; \Longleftrightarrow \; 
\sigma \wedge \sigma^\prime  \! \in \, \rgeo \sigma_j, \sigma \rgeo \; \Longleftrightarrow\;  \sigma^\prime \! 
\in \! T_j^b. $$
This proves that 
\begin{equation} 
\label{aball}
\Gamma (\sigma, r)= T(a) \cap T^b_j \; .
\end{equation}
Conversely, let $j\in \cJ_b$ be such that $h(T^b_j):= 
\max \big\{ d(\sigma_j, \gamma) ;  \gamma \in T_j^b \big\} \geq r/2$. Let $\sigma \in T(a) \cap T^b_j$; then the previous arguments imply (\ref{aball}). Since $T$ is compact, the set $\{ j\in \cJ_b: h(T^b_j) \geq r/2\}$ is finite, which completes the proof of $(i)$.

Let us prove $(ii)$: let $\sigma \in T(a)$, let $r\in (0, 2a)$ and set 
$\delta= {\rm diam} (\Gamma (\sigma, r))$. Then (\ref{aball}) implies that $\Gamma (\sigma, r)$ is compact and there are $\sigma_1, \sigma_2 \in \Gamma (\sigma, r)$ such that $d(\sigma_1, \sigma_2)= \delta$. Observe that it implies  
$$\Gamma (\sigma, r)= \big\{ \sigma^\prime \in T(a): \sigma_1 \wedge \sigma_2 \in \lgeo \rho , \sigma^\prime \rgeo \big\}\; .$$ 
Thus, $\Gamma (\sigma, r)= \overline{\Gamma} (\sigma, \delta)$, that is the closure of $\Gamma (\sigma, \delta)$, and it implies (\ref{diamboule}). The set of all $T(a)$-balls is therefore 
$\bigcup_{q\in \bbQ \cap [0, \infty)} \ccB_{a, q}$, which is a countable set.  
%

Let us prove $(iii)$: $r^\prime<r$ and $\sigma, \sigma^\prime \in T(a)$ and suppose that 
$\Gamma(\sigma^\prime, r^\prime) \cap  \Gamma(\sigma, r)\neq \emptyset$. Then $(i)$ and $(ii)$ implies that $\Gamma (\sigma, r) = \Gamma(\sigma^\prime, r)$, which implies that 
$\Gamma(\sigma^\prime, r^\prime) \subset   \Gamma(\sigma, r)$.
\cqfd

\subsection{Comparison lemmas for Hausdorff measures on real trees.}

Let $(T,d, \rho)$ be a compact real tree. We briefly recall the definition of Hausdorff measures on $T$ and we state two comparison lemmas that are used in the proofs. Let $r_0 \in (0, \infty)$ and 
let $g: [0, r_0) \rightarrow [0, \infty)$ be a function that is assumed to be increasing, continuous and such that $g(0)=0$. For all $\varepsilon \in (0, r_0)$ and all $A \subset T$, we set 
$$ \cH^{(\varepsilon)}_g (A)  = \inf\left\{\sum\limits_{n\in\bbN} g \left({\rm diam}(E_n)\right) \, ; \, A\subset \bigcup\limits_{n\in\bbN} E_n, \,  {\rm diam}(E_n)<\veps\right\} $$
and 
$$ \cH_g (A)= \lim\limits_{\veps \downarrow 0} \uparrow \cH^{(\varepsilon)}_g (A) \; .$$
Under our assumptions, $\cH_g$ is a Borel-regular outer measure : this is the $g$-Hausdorff measure on $T$ (see Rogers \cite{Rogersbook}). The following comparison lemma was first stated for Euclidean spaces by Rogers and Taylor \cite{RoTa61}. 
The proof can be easily adapted to general metric spaces (see Edgar \cite{Edg07}). We include a brief proof of it in order to make the paper self-contained. 
\begin{lem}\label{Comparison1}
Let $(T,d, \rho)$ be a compact rooted real tree.  
Let $\mu$ be a Borel measure on $T$. Let $A$ be a Borel subset of $T$ and let $c\in (0, \infty)$. 
Assume that 
$$\forall\sigma\in A \quad  \limsup\limits_{r\to 0}\frac{\mu\left(B(\sigma,r)\right)}{g(r)}<c  \; .$$
Then, $\mu(A)\leq c \cH_g(A)$. 
\end{lem}

\noindent{\bf Proof.} For any $\veps \! \in \! (0, r_0)$, set 
$$A_\veps= \big\{\sigma  \! \in \! A :  \sup_{ r\in(0,\veps)} \frac{\mu(B(\sigma ,r))}{g(r)}  \! <\!  c \big\} \; .$$ 
Observe that for all $\veps^\prime \! <\!  \veps$, $A_\veps \subset A_{\veps^\prime} \subset A$ and 
$A \! =\! \bigcup_{\veps \in (0, r_0)} 
A_\veps $. Let $(E_n)_{n\in \bbN}$ be a $\veps$-covering of $A_\veps$: namely 
$A_\veps \subset \bigcup_{n \in \bbN}E_n$ and ${\rm diam}(E_n)<\veps$, for all $n\in \bbN$. 
Set $I= \{n\in \bbN: E_n \cap A_\veps\neq \emptyset \}$ and for all $n\in I$, 
fix $\sigma_n \in E_n \cap A_\veps$. Since $g$ is continuous, for all $n\in I$ there exists $r_n \in ( {\rm diam} (E_n), \veps)$ such that 
$$ E_n \subset B(\sigma_n, r_n) \quad \textrm{and} \quad g(r_n) \leq 2^{-n-1} \veps + g( {\rm diam} (E_n)) \; .$$ 
Observe that $\mu (B(\sigma_n ,r_n)) \!<\! cg(r_n)$ and that $A_\veps \! \subset \! \bigcup_{n\in I} B(\sigma_n, r_n) $. Thus, 
\begin{eqnarray*}
\mu(A_\veps) & \leq & \mu \Big( \bigcup_{n\in I}B(\sigma_n ,r_n) \Big)  \leq \sum_{n\in I} \mu (B(\sigma_n ,r_n)) \\
& \leq & \sum_{n\in I} c\, g(r_n)  \leq c\varepsilon +  \sum_{n\in \bbN} c\, g( {\rm diam} (E_n)) \; .
\end{eqnarray*}
Taking the infimum  over all the possible $\veps$-coverings of $A_\veps$ yields 
$$\mu(A_\veps) \leq c\veps + c\cH_g^{(\veps)} (A_\veps) \leq c\veps + c\cH_g (A_\veps) \leq c\veps+ \cH_g(A)\; , $$
which implies the desired result since 
$\mu(A) = \lim_{\veps \downarrow 0} \uparrow \mu(A_\veps)$. 
\cqfd

\medskip

In the next comparison lemma, that seems new to us, we restrict our attention to the level sets of real trees. A more general variant of this result involves a multiplicative constant depending on the gauge function. It has been first stated in Euclidian spaces  by Rogers and Taylor \cite{RoTa61}  (see also Perkins \cite{Per88}) and in general metric spaces (see Edgar \cite{Edg07}). 
\begin{lem}\label{Comparison2}
Let $(T,d, \rho)$ be a compact rooted real tree. Let $a\in (0, \infty)$ be such that the $a$-level set $T(a)$ is not empty. Let $\mu$ be a finite Borel measure on $T$ such that $\mu(T\backslash T(a))= 0$. Let 
$A\subset T(a) $ be a Borel subset and let $c\in (0, \infty)$. 
Assume that 
$$\forall\sigma\in A \quad  \limsup\limits_{r\to 0}\frac{\mu\left(B(\sigma,r)\right)}{g(r)}>c  \; .$$
Then, $\mu(A)\geq  c \cH_g(A)$. 
\end{lem}

\noindent{\bf Proof.} Let $\veps\in(0,(2a)\wedge r_0)$. Let $U$ be an open set of $T$ such that $A\subset U$. 
For all $\sigma\in A$, there exists $r_\sigma\in(0,\veps)$ such that 
$$ \mu\left(\Gamma(\sigma,r_\sigma)\right)=\mu\left(B(\sigma,r_\sigma)\right)>cg(r_\sigma) \quad \textrm{and} \quad \Gamma (\sigma, r_\sigma) \subset U \; .$$ 
Thus, $A\subset \bigcup_{\sigma \in A} \Gamma (\sigma, r_\sigma) \subset U$. Then, 
Proposition \ref{propGamm} $(ii)$ asserts that the set of all $T(a)$-balls is countable and Proposition \ref{propGamm} $(iii)$ asserts that two $T(a)$-balls are either contained one in the other or disjoint. Therefore, there exists $I \subset \mathbb{N}$ and $\sigma_n \in A$, $n\in I$, such that the $\Gamma (\sigma_n, r_{\sigma_n})$, $n \in I$, are pairwise disjoint and $A \subset \bigcup_{n\in I} \Gamma (\sigma_n, r_{\sigma_n}) \subset U$. 
Moreover, by Proposition \ref{propGamm} $(ii)$, ${\rm diam} (\Gamma (\sigma_n, r_{\sigma_n})) \leq r_{\sigma_n}$. Thus, we get 
\begin{eqnarray*}
c\, \cH_g^{(\veps)} (A) & \leq & \sum_{n\in I} c\,  g\big( {\rm diam}(\Gamma(\sigma_n 
,r_{\sigma_n})) \big) \leq \sum_{n\in I}  c\,  g(r_{\sigma_n}) \\
& \leq &  \sum_{n\in I}
\mu \big( \Gamma(\sigma_n ,r_{\sigma_n}) \big) =
\mu \Big( \bigcup_{n\in I} \Gamma(\sigma_n , r_{\sigma_n}) \Big) \leq \mu (U) \; .
\end{eqnarray*}
As $\veps \rightarrow 0$, it entails $c \, \cH_g(A) \leq \mu(U)$, for all open set $U$ containing $A$. Since $\mu$ is a finite Borel measure, it is outer-regular for the open subsets, which implies the desired result. \cqfd

\section{Preliminary results on the Brownian tree.}\label{secBtree}

\subsection{Basic facts on the Brownian excursion.}\label{basicfactsBrownian}

We work under the  excursion measure $\N$ defined in the introduction and $e$ denote the canonical excursion whose duration is denoted by $\zeta$ (see (\ref{loizeta})). We shall denote by $(\cT, d, \rho)$ the compact rooted real tree coded by $e$.

\paragraph{The branching property.} Fix $b\in (0, \infty)$. We discuss here a decomposition of 
$e$ in terms of its excursions above level $b$; this yields a decomposition of the Brownian tree called the branching property. To that end we first introduce the following time change: for all $t\in [0, \infty)$, we set 
\begin{equation}\label{defGb}
 \tau_b(t)= \inf \Big\{ s\in [0, \infty) : \int_0^s\!\!\! \un_{\{e_u\leq b\}} \ddr u >t 
\Big\} \quad \textrm{and} \quad \widetilde{e}_b(t) = e (\tau_b(t)).
\end{equation}
Note that $(\widetilde{e}_b(t))_{t\in [0, \infty)}$ codes the tree below $b$ namely $\{ \sigma\in \cT: d(\rho, \sigma) \leq b\}$ that is the closed ball with center $\rho$ and radius $b$. We denote by $\cG_b$, the 
sigma-field generated by $(\widetilde{e}_b(t))_{t\in [0, \infty)}$ and completed with the $\N$-negligible sets. The approximation (\ref{loctimeexc}) implies that $L^b_\zeta$ is $\cG_b$-measurable.
Then denote by $(\alpha_j,\beta_j)$, $j\in\cJ_b$, the connected components of the time-set $\{s\in [0, \infty) : e(s)>b\}$. Namely, 
$$ \bigcup_{j\in \cJ_b} (\alpha_j, \beta_j) = \{s\in [0, \infty) : e(s)>b\} \; ,$$
and we call  $(\alpha_j, \beta_j)$ the excursion intervals of $e$ above level $b$. 
For all $j\in \cJ_b$, we next set 
$$ l_j^b= L^b_{\tau_{b} (\alpha_j) } \quad \textrm{and} \quad 
\forall s \in [0, \infty), \quad e_j^b (s)= e_{(\alpha_j+s)\wedge\beta_j}-b\; .$$
Then, the $(e_j^b)_{j\in \cJ_b}$ are the excursions of $e$ above level $b$. Recall from (\ref{rephT}) and (\ref{notationNa}) the notation $\N_b= \N (\, \cdot \, | \, \sup e >b)$, that is a probability measure. The \textit{branching property} asserts the following: under $\N_b$ and conditionally on $\cG_b$, the measure
\begin{equation}
\label{branchpropexc}
\cM_b(\ddr l,\ddr e)= \sum_{j\in \cJ_b} \delta_{(l^b_j, e^b_j)}
\end{equation}
is a Poisson point measure on $[0, L^b_\zeta] \! \times \! \bC^0$ 
with intensity $\un_{[0, L^b_\zeta]} (l) \ddr l \,  \N (\ddr e) $. 

\medskip

The following decomposition of $e$ is interpreted in terms of the Brownian tree $\cT$ as follows. 
Recall that $p: [0, \zeta] \rightarrow \cT$ stands for the canonical projection. Then for all $j\in \cJ_b$, we set 
$$ \sigma_j=p(\alpha_j)= p(\beta_j)\, , \quad \cT^{o,b}_j= p\big( (\alpha_j, \beta_j)\big) \quad \textrm{and} \quad \cT_j^b = p\big( [\alpha_j, \beta_j] \big)\; .$$
Then, we easily check that the $\cT^{o,b}_j$, $j\in \cJ_b$, are the connected components of the open subset $\{ \sigma \!\in \!\cT: d(\rho, \sigma) \!>\! b\}$ and that $\cT^{o,b}_j= \cT^b_j \backslash \{ \sigma_j\}$. 
Namely, the $(\cT_j, d, \sigma_j)$, $j\in \cJ_b$ are the subtrees above level $b$ of $\cT$ as 
introduced in Section \ref{subsecballs}. Moreover note that for all $j\in \cJ_b$, 
the rooted compact real tree $(\cT_j, d, \sigma_j)$ is isometric to the tree coded by the excursion $e_j^b$. We next use this and Proposition 
\ref{propGamm} to discuss the balls in a fixed level of $\cT$. 

   To that end, we fix $a, r\in (0, \infty)$ such that $a>r/2$ and we conveniently set 
$b \!=\! a\!-\!r/2$. Recall that $\cT(a)= \{ \sigma \! \in \! \cT: d(\rho, \sigma) \!=\! a\}$ and that for all $\sigma \! \in \! \cT(a)$, we have set $\Gamma (\sigma, r) \! =\!  \cT(a) \cap B(\sigma, r)$ that is the ball in $\cT(a)$ with center $\sigma$ and radius $r$. We also recall that $\mathscr{B}_{a,r}= \{ \Gamma (\sigma, r); \sigma \! \in \! \cT(a)\}$ stands for the set of all $\cT(a)$-balls with radius $r$. By Proposition \ref{propGamm}, $\mathscr{B}_{a,r}$ is a finite set and that 
$$ \mathscr{B}_{a,r}= \big\{ \cT(a) \cap \cT^b_j \, ; \; j\! \in \! \cJ_b : h(\cT^b_j) \geq r/2  \big\} \; , $$
where the trees $(\cT^b_j, d, \sigma_j)$, $j\! \in \! \cJ_b$, are the subtrees of $\cT$ above level $b$ as previously defined; here $h(\cT^b_j)= \sup_{\sigma \in \cT^b_j} d(\sigma_j, \sigma) $ stands for 
the total height of $\cT^b_j$. Note that $ h(\cT^b_j)= \sup e^b_j$ that is maximum of the excursion corresponding to $\cT^b_j$, as explained above. 

  Then, we set $Z_{a,r}= \# \mathscr{B}_{a,r}$, that is the number of $\cT(a)$-ball with radius $r$. 
 Assume that $Z_{a,r} \geq 1$. We then define the indices $j_1, \ldots, j_{Z_{a,r}} \in \cJ_b$ by 
$$  \{ j_1, \ldots, j_{Z_{a,r}} \} = \big\{ j \! \in \! \cJ_b : h(\cT^b_j) \geq r/2 \big\} \quad \textrm{and} \quad \alpha_{j_1}< \ldots < \alpha_{j_{Z_{a,r}}} \; .$$ 
and we set 
\begin{equation}\label{enumTaballs}
 \forall i\in \{ 1, \ldots , Z_{a,r} \}, \quad \Gamma_i := \cT(a) \cap \cT^b_{j_i} \; .
 \end{equation}
Namely $\mathscr{B}_{a, r} = \big\{ \Gamma_i \, ; \; 1\leq i\leq Z_{a,r}\big\}$ is the set of the $\cT(a)$-balls with radius $r$ listed in their order of visit by the excursion $e$ coding $\cT$.

\begin{lem}
\label{Tabrownian} Let $a, r\in (0, \infty)$ such that $a>r/2$. Let $\big\{ \Gamma_i \, ; \; 1\leq i\leq Z_{a,r}\big\}$ is the set of the $\cT(a)$-balls with radius $r$ listed in their order of visit as explained above. Then the following holds true. 
\begin{itemize}
\item[(i)] Under $\N_a= \N(\, \cdot \, | \, \sup e >a)$, $Z_{a,r}$ has a geometric law with parameter $2a/r$. Namely, 
$$ \forall k\geq 1, \quad \N_a [Z_{a,r}= k]=  \big( 1\! -\! \frac{_r}{^{2a}}\big)^{k-1}\frac{_r}{^{2a}} \; .$$
\item[(ii)] For all $k\geq 1$, under $\N_a (\, \cdot \, | \, Z_{a,r}\! = \! k)$, the r.v.~$(\ell^a(\Gamma_i))_{1\leq i\leq k}$ are independent and exponentially distributed with mean $r/2$. 
\end{itemize}
\end{lem}

\noindent{\bf Proof.}
Let $a\in(0,\infty)$ and denote $b=a\!-\!r/2$. Let $k\geq 1$ and $F_1, \ldots F_k: \bC^0 \rightarrow [0, \infty)$ be measurable functionals. Recall from (\ref{rephT}) that $\N (\sup e \geq r/2)= 2/r$. Then, the definition of the $j_i$ combined with 
the branching property and basic results on Poisson point measures entail
\begin{equation}\label{branchingindep}
\N_b \Big[  \un_{\{Z_{a,r} =k\}} \!\! \prod_{1\leq i\leq k}\!\! F_i (e_{j_i}^b) \,  \Big| \, \cG_b\Big]= \frac{(\frac{2}{r} L^b_\zeta)^k}{k!} e^{-\frac{2}{r} L^b_\zeta}  \prod_{1\leq i\leq k} \N_{r/2} \big[F_i (e) \big] \; .
\end{equation}

Then recall  (\ref{Laplaceella}) that implies that  $L^b_\zeta$ under $\N_b$ is exponentially distributed with mean 
$b$. Thus, 
$$ \frac{1}{{k!}} \N_b \big[ \,  \big(\frac{_2}{^r} L^b_\zeta \big)^k  e^{-\frac{2}{r} L^b_\zeta} \, \big]=  
\frac{(\frac{2}{r} b)^k}{(1+ \frac{2}{r} b )^{k+1}}= \frac{_r}{^{2a}} \big(1- \frac{_r}{^{2a}} \big)^k \; ,$$
because $b= a\!-\!r/2$ and  $(1+ \frac{2}{r} b)^{-1}= r/(2a)$. It implies 
$$ \N_b \Big[  \un_{\{Z_{a,r} =k\}} \!\! \prod_{1\leq i\leq k}\!\! F_i (e_{j_i}^b)\Big]= \frac{_r}{^{2a}} \big(1- \frac{_r}{^{2a}} \big)^k  \prod_{1\leq i\leq k} \N_{r/2} \big[F_i (e) \big] \; .$$
Next observe that $\N_b$-a.s.~$\un_{\{ \sup e >a\}}= \un_{\{ Z_{a,r} \geq 1\}}$. Thus, we get 
\begin{eqnarray}
\label{poipoifor}
 \N_a \Big[ \un_{\{Z_{a,r} =k\}} \!\! \prod_{1\leq i\leq k}\!\! F_i (e_{j_i}^b)\Big] &=& \frac{_a}{^b}  \N_b \Big[  \un_{\{Z_{a,r} =k\}} \!\! \prod_{1\leq i\leq k}\!\! F_i (e_{j_i}^b)\Big] \nonumber \\
 & =& \frac{_r}{^{2a}} \big(1- \frac{_r}{^{2a}} \big)^{k-1}  \prod_{1\leq i\leq k} \N_{r/2} \big(F_i (e) \big)
\end{eqnarray}
because $a/b\! = \! (1\! -\! \frac{r}{2a})^{-1}$. Recall that (\ref{Laplaceella}) implies that under $\N_{r/2}$, $\ell^{r/2} (\cT) \! =\! 
 L_\zeta^{r/2} $ is exponentially distributed with mean $r/2$. By taking 
$F_i(e) \! =\!  f_i (L^{r/2}_\zeta)$ in (\ref{poipoifor}) we then get 
$$  \N_a \Big[  \un_{\{Z_{a,r} =k\}} \!\! \prod_{1\leq i\leq k}\!\! f_i \big(
\ell^a(\Gamma_i) \big) \Big]= \frac{_r}{^{2a}} \big(1- \frac{_r}{^{2a}} \big)^{k-1} \!\!\!  \prod_{1\leq i\leq k} \int_0^\infty\!\!\!  f_i(s)\frac{_2}{^r}e^{-\frac{_2}{^r}s} \ddr s \; , $$
with entails the desired result. 
\cqfd

\paragraph{Ray-Knight theorem under $\N$.} We first recall the definition of Feller diffusion, namely a 
Continuous States space Branching Process (CSBP) with branching mechanism 
$\psi (\lambda)= \lambda^2$. Let $x\in [0, \infty)$ and let $(Y^x_a)_{a\in [0, \infty)} $ be a $[0, \infty)$-valued continuous process defined on the probability space 
$(\Omega, \cF, \P)$. It is a Feller diffusion with branching mechanism $\psi (\lambda)= \lambda^2$ and initial value $Y^x_0= x$ 
if it is a Markov process such that 
$$ \E \big[ \exp (-\lambda Y^x_{a+a^\prime})\,  \big| Y^x_a \big]= \exp \Big( \!-\! \frac{Y^x_a \lambda }{1+ a^\prime\lambda } \Big) \; , \quad a, a^\prime, \lambda \in [0, \infty)\; .$$
Recall notation $\N_a= \N (\cdot \big| \, \sup e >a)$ and
$\cG_a$ for the sigma-field generated by the excursion $\widetilde{e}_a$ defined in (\ref{defGb}). Recall that $\ell^a(\cT)= L^a_\zeta$, the total mass of the local-time measure at level $a$, is $\cG_a$-measurable.

\noi
We shall use the following statement of Ray-Knight theorem. Let $a\in (0, \infty)$. 
\begin{itemize}
\item[$(i)$] $\N_a [\exp(-\lambda \ell^a (\cT)) ]= \frac{1}{1+ a\lambda}$. 

\item[$(ii)$] Under $\N_a$ and conditionally given $\cG_a$, the process $(\ell^{a+a^\prime} (\cT))_{a^\prime\in [0, \infty)}$ is a Feller diffusion with branching mechanism $\psi(\lambda)= \lambda^2$ and initial value $\ell^a(\cT)$. 
\end{itemize}
This is an immediate consequence of the Ray-Knight theorem for standard Brownian motion and of the Markov property under $\N$ : see \cite{Blubook} III 3 and VI 2.10. 

Combined with the branching property, the above Ray-Knight theorem,  has the following consequence. Let us recall that we  enumerate the $\cT(a)$-balls of $\ccB_{a,r}$ as $\{\Gamma_i,  1\leq i \leq Z_{a,r}\}$ (see (\ref{enumTaballs})).  Let $\Gamma$ such a $\cT(a)$-ball. For $a^\prime\geq 0$, we define
\begin{equation}\label{postboule}
\Gamma^{a+a^\prime}=\left\{\sigma\in\cT(a+a^\prime) \  \exists \sigma^\prime\in\Gamma : \sigma^\prime\in\llbracket\rho,\sigma\rrbracket\right\},
\end{equation}
the set of vertices at level $a+a^\prime$ that have an ancestor in $\Gamma$ (notice that $\Gamma^a=\Gamma$). The  following lemma is a straightforward consequence of  Ray-Knight theorem.
\begin{lem}\label{branchingRK}
Let $a\in(0,\infty)$, $r\in[0,2a]$. Let $\{\Gamma_i,   1\leq i\leq Z_{a,r}\}$ the set of $\cT(a)$-balls of radius~$r$. Under $\N_a$ conditionally on $\cG_a$, the processes $\left(\ell^{a+a^\prime}(\Gamma_i^{a+a^\prime}), a^\prime\geq 0\right),  1\leq i\leq Z_{a,r},$ are independent Feller diffusions started at $\left(\la(\Gamma_i)\right), 1\leq i \leq Z_{a,r}$.
\end{lem}
\noindent{\bf Proof.}
Recalling for $b=a\!-\!r/2$ the decomposition (\ref{enumTaballs}), we see that 
\begin{equation}\label{enumpostboules}
 \forall i\in \{ 1, \ldots , Z_{a,r} \}, \quad \Gamma_i^{a+a^\prime} := \cT(a+a^\prime) \cap \cT^b_{j_i} \; .
\end{equation}
Hence, one can use (\ref{branchingindep}), and the Ray-Knight theorem (see $(ii)$ above) to get the desired result.~\cqfd

\paragraph{Spinal decomposition.} We recall another decomposition of the Brownian tree called \textit{spinal decomposition}. This is a consequence of Bismut's decomposition 
of the Brownian excursion that we recall here.  

Let $X$ be a real valued process defined on $(\Omega, \cF, \P)$ such that $(\frac{1}{\sqrt{2}} X_t)_{t\in [0, \infty)}$ is distributed as a standard Brownian motion with initial value $0$. Let $X^\prime$ be an independent copy of $X$ on  $(\Omega, \cF, \P)$. We fix $a\in (0, \infty)$ and we set 
$$ T_a= \inf \{ t\in [0, \infty): X_t= -a \}\quad  \textrm{and} \quad T^\prime_a= \inf \{ t\in [0, \infty): X^\prime_t= -a \}\; .$$
We next set for any $s\in [0, \infty)$, 
$$ \check{e}^t_s = e_{(t-s)_+}\quad \textrm{and} \quad \hat{e}^t_s=e_{t+s} \; . $$
Then the Bismut's identity (see \cite{Bis85} or \cite{LG93}) states that for any non-negative measurable functional $F$ on $(\bC^0)^2$,
\begin{equation}
\label{Bibifoc}
 \N \Big[ \int_0^\zeta \!\!\! \ddr L^a_t \, F\big(  \check{e}^t \, ; \, \hat{e}^t \big)\Big]= \E \big[ F( a+X_{\cdot \wedge T_a} ; a+X^\prime_{\cdot \wedge T^\prime_a })\big] \; .
\end{equation} 
We derive from (\ref{Bibifoc}) an identity involving the excursions above the infimum of $\hat{e}^t$ and $\check{e}^t$. To that end, we introduce the following.  
Let $h: [0, \infty ) \rightarrow [0, \infty)$ with compact support. We define a point measure point $\cN (h)$  as follows: set $\underline{h} (t)= \inf_{[0, t]} h$ and denote by $(g_i, d_i)$, $i\in \cI(h)$ the excursion intervals of $h-\underline{h} $ away from $0$ that are the connected component of the open set $\{ t \geq 0: h(t)-\underline{h} (t) >0 \}$. For any $i \in \cI (h)$, set 
$h^i (s)  = ((h- \underline{h}) ( (g_i +s) \wedge d_i) \, , \, s \geq 0)$. We then define $\cN(h)$ as the point measure on $[0, \infty) \times \bC^0$ given by 
$$ \cN (h)= \sum_{ i\in \cI (h)} \delta_{(h(g_i) , h^i) } \; .$$
Then, for any $t, a\in (0, \infty)$, 
\begin{equation}
\label{spinaldef}
\cN_{t} := \cN ( \check{e}^t ) + \cN(\hat{e}^t) =: \sum_{j\in \cJ_t } \delta_{(h^t_j, e^{t,j} )} 
\end{equation}
and 
\begin{equation}
\label{Nstaradef}
 \cN_{a}^* := \cN ( a+X_{\cdot \wedge T_a}) + \cN(a+X^\prime_{\cdot \wedge T^\prime_a }) =: \sum_{j\in \cJ^*_a } \delta_{(h^{*}_j, e^{*j})}\; .
\end{equation} 
We deduce from (\ref{Bibifoc}) that 
for all any $a$ and for all nonnegative measurable function $F$ on the set of positive measures on $[0, \infty) \times \cC^0$, one has 
\begin{equation}      
\label{Bibiexcu}
\N \Big[ \int_{0}^\zeta \!\!\! \ddr L^a_t  \; F \big(  \cN_t  \big)  \Big] =  \E \big[ F (\cN^*_a ) \big] 
\end{equation}
and as consequence of It\^o's decomposition of Brownian motion above its infimum, $\cN^*_a$ is a Poisson point measure on $[0, \infty) \times \bC^0$ with intensity 
$2 \un_{[0, a]} (h) \ddr h \, \N (\ddr e)$.\\

Let us interpret this decomposition in terms of the Brownian tree. Choose $t \in (0, \zeta)$ such that $e_t=a$ and set $\sigma = p(t) \in \cT$ (namely $\sigma \in \cT(a)$). 
Then, the geodesic $\lgeo \rho, \sigma \rgeo$ is interpreted as the ancestral line of $\sigma$. Let us denote by $\cT_{j}^{o}$, $j \in \cJ$, the connected components of the open set $\cT   \backslash \lgeo \rho , \sigma \rgeo$ and denote by $\cT_{j}$ the closure of 
$\cT_j^o$. Then, there exists a point $\sigma_j \in \lgeo \rho , \sigma \rgeo$ such that $\cT_{j} = \{ \sigma_{j} \} \cup \cT_{j}^{o}$. 
Recall notation $(h^{t}_{j} , e^{t,j})$, $j \in \cJ_t$ from (\ref{spinaldef}). The specific coding of $\cT$ by $e$ entails that for any $j \in \cJ$ there exists a unique $j^\prime \in \cJ_{t}$ such that $ d(\rho ,\sigma_{j})= h^{t}_{j^\prime} $ and such that the rooted compact real tree $(\cT_{j}, d, \sigma_{j})$ is isometric to the tree coded by $e^{t,j^\prime}$.

Recall that $p(t)= \sigma$. We fix $r , r^\prime \in [0, 2a)$ such that $r^\prime \leq r$. We now compute the mass of the ring $B(\sigma, r) \setminus B(\sigma , r^\prime)$ in terms of $\cN_t$. 
First, observe that for any $s \in [0,\zeta ]$ such that $e_s = a $, we have 
$$r^\prime \leq d(s, t) < r \Longleftrightarrow a-(r^\prime/2) \geq \inf_{u\in [s\wedge t, s\vee t ]} e_u  >a-(r/2)\; .$$
We then get 
\begin{equation}
\label{linloc}
 \ell^a \big( B(\sigma, r  )\backslash B(\sigma, r^\prime) \big) = \sum_{j \in \cJ_t} \un_{( a-\frac{r}{2} \, , \, a-\frac{r^\prime}{2}]} (h^t_j) \, 
L^{a-h^t_j}_{\zeta^t_j} (t,j)  \; , 
\end{equation}
where $L^{a-h^t_j}_{\zeta^t_j} (t,j) $ stands for the local time at level $a-h^t_j$ of the excursion $e^{t, j}$.

Then, for any $a \in (0, \infty)$ and any $ r, r^\prime \in (0, 2a)$ such that $r^\prime \leq r$, we also set 
\begin{equation}
\label{Lambdaadef}
\Lambda^a_{r^\prime, r}  = \sum_{j \in \cJ^*_a} \un_{(a-\frac{r}{2} \, , \,  a-\frac{r^\prime}{2}]} (h^*_j) \, L_{\zeta^*_j}^{a-h^*_j}  \; , 
\end{equation}  
where, $L_{\zeta^*_j}^{a-h^*_j}$ stands for the local time at level $a-h^*_j$ of the excursion $e^{*  j}$ defined in (\ref{Nstaradef}). 
Then, (\ref{Bibiexcu}) implies that for any $a \in (0, \infty)$ and for all non-negative measurable $F$ function 
\begin{equation}\label{bibiprovisoire}
\N \Big[\int_{\cT} \!\! \ell^a (\ddr \sigma) \, F\big( \,  \ell^a \big( B(\sigma, r  )\backslash B(\sigma, r^\prime)  \big)\,  ; \, 0 \! \leq \! r^\prime \! \leq \! r \! \leq \! 2a \, \big) \Big] =
\E \big[ F\big( \Lambda^a_{r^\prime , r}\, ; \,  0 \! \leq \! r^\prime \! \leq \! r \! \leq \! 2a \, \big)\big]
\end{equation}

On the right-hand-side, the dependency with respect to the level $a$ is a bit artificial. Indeed, for $a\in(0,\infty)$, the Poisson point measure $\cN_a^*(\ddr h\ddr e)$ has its law invariant under the transformation $(h,e)\mapsto(a-h,e)$. Thus, let us consider on $(\Omega,\cF,\P)$ a new Poisson point measure $\cM^*=\sum\limits_{j\in\cI^*}\delta_{(h^*_j,e^*_j)}$ with intensity $2\ddr h\N(\ddr e)$ (we abuse notations and keep the notation $(h^*_j,e^*_j)$ for the atoms).  We  set
\begin{equation}
\label{Lstardef}
\Lambda^*_{r^\prime, r}  = \sum_{j \in \cI^*} \un_{[\frac{r^\prime}{2} \, , \, \frac{r}{2})} (h^*_j) \, L_{\zeta^*_j}^{h^*_j}, 
\end{equation}  
where $L_{\zeta^*_j}^{h^*_j}$ stands for the local time at height $h^*_j$ for the excursion $e^*_j$. One can now rewrite (\ref{bibiprovisoire}) as

\begin{eqnarray}\label{spinalrings}
& &\hspace{-20mm} \N \Big[\int_{\cT} \!\! \ell^a (\ddr \sigma) \, F\big( \,  \ell^a \big( B(\sigma, r  )\backslash B(\sigma, r^\prime)  \big)\,  ; \, 0 \! \leq \! r^\prime \! \leq \! r \! \leq \! 2a \, \big) \Big]  \nonumber  \\
& & \hspace{30mm}=  
\E \big[ F\big( \Lambda^*_{r^\prime , r}\, ; \,  0 \! \leq \! r^\prime \! \leq \! r \! \leq \! 2a \, \big)\big]\; .
\end{eqnarray}

The law of the $\Lambda^*_{r^\prime , r}$ is quite explicit as shown by the following lemma. 

\begin{lem}
\label{otarieimperiale}
Let $0 \! \leq \! r_n \! \leq \! r_{n-1} \! \leq \! \ldots \! \leq \! r_1 \leq 2a$. Then,  
$$ \Lambda_{r_n, r_{n-1} }^* \, , \,  \Lambda_{r_{n-1}, r_{n-2} }^*  \, , \,  \ldots \, , \,  \Lambda_{r_2, r_{1} }^* $$
are independent. Moreover, for any $ 0 \! \leq \! r^\prime \! \leq \! r \! \leq \! 2a$, 
$$ \forall y\in (0, \infty) \quad  \P \big( \Lambda^*_{r^\prime, r}  >y \big)= \left( 1-\frac{{r^\prime}}{r} \right)^{\! \! 2}  \frac{2y}{r} \, e^{-2y/r} +  \left( 1-\Big(\frac{r^{\prime}}{r} \Big)^2 \right) 
e^{-2y/r}, $$
and $\P (\Lambda^*_{r^\prime, r} =0)= (r^\prime/r)^2$. 
\end{lem}

\noindent{\bf Proof.}
The intervals $[r_{k+1}/2,r_{k}/2)$ being pairwise disjoint, the independence of the increments is a straightforward consequence of the properties the  Poisson point measure $\cM^*$. Using Campbell formula and (\ref{Laplaceella}), we compute, for all $\lambda\geq0$,
\begin{eqnarray*}
\E\left[e^{-\lambda \Lambda^*_{2r^\prime,2r}}\right] &=& \exp\left(-\int_{r'}^r 2\ddr h \N\left[1-e^{-\lambda \ell^{h}(\cT)}\right]\right) \\
&= &\exp\left(-\int_{r'}^r 2\ddr h\frac{\lambda}{1+h\lambda}\right)=\left(\frac{1+r'\lambda}{1+r\lambda}\right)^2.
\end{eqnarray*}
Thus, $\Lambda^*_{2r^\prime,2r}\overset{(law)}{=}X_1+X_2$, where $X_1$ and $X_2$ are i.i.d random variables where $$\E\left[e^{-\lambda X_1}\right]=\frac{r'}{r}+\left(1-\frac{r'}{r}\right)\frac{1}{1+r\lambda}.$$
Thus, $X_1=0$ with probability $r'/r$ and conditionally on being non-zero, it is exponentially distributed with mean $r$. Thus, for $y>0$, 
\begin{align*}
\P\left(\Lambda^*_{2r^\prime, 2r}>y\right)&=2\P\left(X_1=0 ;\ X_2>y\right)+\P\left(X_1>0\ ;\ X_2>0\ ;\ X_1+X_2>y\right)\\
&=2\frac{r'}{r}\P(X_1>y)+\left(1-\frac{r'}{r}\right)^2\P\left(Z>y\right),
\end{align*}
where $Z$ has law ${\rm Gamma}(2,1/r)$. The result proceeds now from 
elementary computations. 
\cqfd

\subsection{Estimates.}\label{secestim}
The following elementary computation is needed twice in our proofs.
\begin{lem}\label{moment4}
Let $(X_n)_{n\geq 1}$ a sequence of i.i.d real valued random variables on $(\Omega,\cF,\P)$, with mean $0$ and a moment of order $4$. Let $Z$ be a random variable taking its values in $\mathbb{N}$, independent of the sequence $(X_n)$. Then
$$\E\left[ \left(X_1+ X_2 + \cdots + X_Z\right)^4\right]\leq 3\E[X_1^4]\E\left[Z^2\right].$$ 
Moreover, the following holds : $\E\left[\left(X_1-\E[X_1]\right)^4\right]\leq 2\E[X_1^4]$.
\end{lem}
\noindent{\bf Proof.} One has 
$$\E\left[ \left(X_1+ X_2 + \cdots + X_Z\right)^4\mid Z\right]=
\sum\limits_{1\leq\!i_1,i_2,i_3,i_4\leq\!Z}
\E\left[X_{i_1}X_{i_2}X_{i_3}X_{i_4}\right].$$
When $(i_1,i_2,i_3,i_4)$ contains an index that is distinct of the three others, then the contribution of the corresponding term will be null. Thus the latter mean equals $Z\E[X_1^4]+3Z(Z\!-\!1)\E[X_1^2]^2\leq3 Z^2\E[X_1^4]$ (using Jensen's inequality). The second statement follows from 
\begin{eqnarray*}
\E\left[\left(X_1-\E[X_1]\right)^4\right] &= &
\E\left[\left(X_1-\E[X_1]\right)^4\un_{\{X_1\geq \E[X_1]\}}\right] + \E\left[\left(X_1-\E[X_1]\right)^4\un_{\{X_1< \E[X_1]\}}\right] \\
& \leq & \E\left[X_1^4\right] +  \E[X_1]^4,
\end{eqnarray*}
and using Jensen's inequality. 
\cqfd

\bigskip


We explained in Section \ref{basicfactsBrownian} the link between the process $(\ell^a(\cT), a\in(0,\infty))$ and the  Feller diffusion, for which we provide here some basic estimates.

\begin{lem}\label{GDFeller}
Let $(Y^x_a)_{a\geq 0}$ be a Feller diffusion starting at $x\geq 0$, defined on $(\Omega,\cF,\P)$.
For all $x,y\in[0,\infty)$, for all $a\in(0,\infty)$, the following inequalities hold :
\begin{enumerate}
\item[(i)] If $y\leq x$, then   $\P\left(\inf\limits_{b\in[0,a]}Y_b^x  \leq y\right)\leq \exp\left(-\frac{1}{a}(\sqrt{x}-\sqrt{y})^2\right)$.

\item[(ii)] If $y\geq x$, then  $\P\left(\sup\limits_{b\in[0,a]}Y_b^x \geq y\right)\leq \exp\left(-\frac{1}{a}(\sqrt{y}-\sqrt{x})^2\right).$
\end{enumerate} 
\end{lem}
\noindent{\bf Proof.}
Let us prove $(i)$. Recall that for all $x,b,\lambda\in[0,\infty)$, $\E\left[e^{-\lambda Y^x_b}\right]=\exp\left(-\frac{\lambda x}{1+b\lambda}\right)$. Thus, for fixed $a\in(0,\infty)$,  and for  $\lambda\in[0,\frac{1}{a})$, we set
\begin{equation}
\forall b\in[0,a], \quad M_b^{(\lambda,x)}:=\exp\left(\!-\frac{\lambda Y^x_b}{1-b\lambda}\right).
\end{equation}
We stress that for $b\in[0,a]$, one has $1-b\lambda\geq 1-a\lambda>0$, and one can compute 
$$\E[ M_b^{(\lambda,x)}]=\exp\left(-\frac{\lambda}{1-b\lambda}x\Big{/}\left(1+\frac{b\lambda}{1-b\lambda}\right)\right)=e^{-\lambda x}. $$
Combined with the Markov property, this entails that $(M_b^{(\lambda,x)}, b\in[0,a])$ is a martingale. Moreover, on $\{\inf\limits_{b\in[0,a]}Y_b^x  \leq y\}$, one has $\inf\limits_{b\in[0,a]}\frac{\lambda Y_b^x}{1-b\lambda}\leq \inf\limits_{b\in[0,a]}\frac{\lambda Y_b^x}{1-a\lambda}\leq\frac{\lambda y}{1-a\lambda}$. Hence, the maximal inequality for sub-martingales entails
\begin{eqnarray*}
\P\left(\inf\limits_{b\in[0,a]}Y^x_b \leq y\right)& \leq & \P\left(\sup\limits_{b\in[0,a]}
M_b^{(\lambda,x)}\geq e^{-\frac{\lambda y}{1-a\lambda}}\right) \\
&\leq & e^{\frac{\lambda y}{1-a\lambda}}\E\left[ M_a^{(\lambda,x)}\right]=\exp\left(\frac{\lambda y}{1-a\lambda}-\lambda x\right).
\end{eqnarray*}
The reader can check using elementary computations that the function $\lambda\mapsto\frac{\lambda y}{1-a\lambda}-\lambda x$ has a negative minimum on $(0,1/a)$ at the value $\lambda=\frac{1}{a}\left(1-\sqrt{\frac{y}{x}}\right)$, and this minimum is $-\frac{1}{a}(\sqrt{x}-\sqrt{y})^2$, which completes the proof.\\

In order to prove $(ii)$, one could extend the definition of $(M_b^{(\lambda,x)}, b\in[0,a])$ for $\lambda\in(-1/a,0)$. In what follows, we use a simpler argument. Let us begin with the following remark: let  $b\in(0,\infty)$, let $\mathcal{E}$ be a r.v.~on $(\Omega,\cF,\P)$ that is exponentially distributed with mean $b$,  then for all $\lambda\geq 0$, $\E[e^{-\lambda \mathcal{E}}]=\frac{1}{1+b\lambda}$, and this Laplace transform remains finite for $\lambda\in(-1/b, 0)$. Moreover, one can plainly check that for $x,b\in(0,\infty)$, $Y_b^x$ has the same law as $\sum\limits_{i=1}^N\mathcal{E}_i$, where the $\mathcal{E}_i$ are independent copies of $\mathcal{E}$ and $N$ is an independent Poisson r.v. with mean $x/b$. Thus, one has
\begin{equation}
\forall \mu\in(0,1/b),\quad \E\left[e^{\mu Y_b^x}\right]=\exp\left(\frac{\mu x}{1-\mu b}\right).
\end{equation}
The Feller diffusion $(Y_b^x, b\geq 0)$ is a martingale, so by convexity $(e^{\mu Y_b^x}, b\geq 0)$ is a submartingale. Thus, for all $\mu\in(0,1/a)$, and $y\geq x\geq 0$, one has
\begin{eqnarray*}
\P\left(\sup\limits_{b\in[0,a]}Y^x_b \geq y\right)& \leq & \P\left(\sup\limits_{b\in[0,a]}
e^{\mu Y_b^x} \geq e^{\mu y}\right) \\
& \leq & e^{-\mu y}\E\left[e^{\mu Y_a^x}\right]=\exp\left(\frac{\mu x}{1-a\mu}-\mu y\right),
\end{eqnarray*}
and the result follows by optimizing  the same function as before. 
\cqfd

\bigskip

\noi The next result is a corollary of Lemma \ref{GDFeller} $(ii)$.

\begin{lem}\label{GDella}
Let $m\in(0,1/2)$. For all $y\in(0,\infty)$, 
$$\N\left(\sup\limits_{b\in[m,m^{-\!1}]}\lb(\cT)>y\right)\leq (2/m)\exp\left(-my/2\right).$$
\end{lem}
\noindent{\bf Proof.}
Let $m\in(0,1/2)$ and recall from (\ref{defGb}) the definition of $\cG_m$. As recalled in Section \ref{basicfactsBrownian},  under  $\N_{m}$, conditionally on $\cG_m$, the process $(\lb(\cT), b\geq m)$ is a Feller diffusion started at $\ell^{m}(\cT)$. Hence, conditioning with respect to  $\cG_m$ and using Lemma \ref{GDFeller} $(ii)$, we get 
$$\N_m\left(\sup\limits_{b\in[m,m^{-\!1}]}\lb(\cT)>y\right)
 \leq \N_m\left[\exp\left(-m\left(\sqrt{y}-\sqrt{\ell^m(\cT)}\right)^2\right)\right].$$
Expanding $(\sqrt{u/2}-\sqrt{2v})^2$, one shows that for all $u,v\geq 0$, $\left(\sqrt{u}-\sqrt{v}\right)^2\geq u/2-v$. Thus, $\N_m\left(\sup\limits_{b\in [m,m^{-\!1}]}\ell^m(\cT)>y\right)\leq \exp\left(-m\frac{y}{2}\right)\N_m\left[e^{m\ell^m(\cT)}\right]$. Recalling from (\ref{Laplaceella}) that under $\N_m$,  $\ell^m(\cT)$ is exponentially distributed with mean $m$, we get $\N_m\left[e^{m\ell^m(\cT)}\right]=(1-m^2)^{\!-1}\leq 2$, because  $m<1/2$. This entails the desired result, recalling that $\N_m(\cdot)=m\N\left(\cdot\un_{\{h(\cT)>m\}}\right)$ and that the events $\{h(\cT)>m\}$ and $\{\ell^m(\cT)>0\}$ are equal, up to a $\N$ negligible set. 
\cqfd

\paragraph{Estimates for small balls.}
We consider here a level $a\in(0,\infty)$ and recall   that $\cT(a)$ is the $a$-level set of the Brownian tree $\cT$. If  $r\in [0,2a]$, we recall from (\ref{Taboule}) the notation $\Gamma(\sigma,r)$ for the $\cT(a)$-ball of radius $r$ and center $\sigma\in\cT(a)$, the set of $\cT(a)$-balls of radius $r$ being denoted $\ccB_{a,r}$.  Let $\Gamma$ be a $\cT(a)$-ball of radius $r^\prime$, where $r^\prime\in[0,2a]$.  From Proposition \ref{propGamm} $(iii)$, we know that if $r\in[r^\prime,2a]$, there exists a unique $\cT(a)$-ball of radius $r$ that contains $\Gamma$, and we shall denote  this "enlarged" ball by
\begin{equation}\label{enlargedball}
\Gamma[r]:=\Upsilon \quad {\rm where\ } \Upsilon\in\ccB_{a,r} {\rm\  and \ } \Gamma\subset\Upsilon.
\end{equation}

We consider positive real numbers $r_1>r_2>\ldots>r_n>0$, and $\veps_1>\ldots>\veps_{n\!-\!1}>0$, where $n\in\mathbb{N}^*$. We set $\br=\left\{r_1,\ldots,r_n\right\}$ and $\beps=\left\{\veps_1,\ldots,\veps_{n\!-\!1}\right\}$.
We shall say that $\Gamma$, a  $\cT(a)$-ball of radius $r_n$,  is $(\br,\beps)$-small if and only if for all  $1\leq k\leq n\!-\!1$, the enlarged ball of radius $k$ has a local time smaller than $\veps_k$, namely
\begin{equation}\label{condballsmall}
\forall k\in\{1,\ldots,  n\!-\!1\}  \quad\la\left(\Gamma[r_k]\right)\leq \veps_k.
\end{equation}
We denote by $S_{a,\br,\beps}$ the  total number of such $(\br,\beps)$-small balls at level $a$: 
\begin{equation}\label{defZs}
S_{a,\br,\beps}:=\sum\limits_{\Gamma\in\ccB_{a,r_n}}\un_{\{\Gamma {\rm\  is \ }(\br,\beps)-{\rm small} \}}.
\end{equation}
To control that number, we introduce
\begin{equation}\label{defmu}
\mu(\br,\beps):=\N\left[S_{r_1/2,\br,\beps}\right].
\end{equation}
Let us stress that its definition does not depend on $a$.

\begin{lem}\label{Z4}
Let $a\in(0,\infty)$, $\br=\{r_1,\ldots,r_n\}$, and $\beps=\left\{\veps_1,\ldots,\veps_{n\!-\!1}\right\}$, where $r_1>\ldots>r_n>0$, and $\veps_1>\ldots>\veps_{n\!-\!1}>0$. There exists a constant $c_0\in(0,10^4]$ such that  if $a/r_1>1$ and $r_1/r_n>2$, 

$$\N\left[\left(S_{a,\br,\beps}-\mu(\br,\beps)\la(\cT)\right)^{\!4}\right]\leq c_0 a\ \frac{r_1^2}{r_n^4}.$$
\end{lem}

\noindent{\bf Proof.}
Let $a,\br,\beps$ as above. From Proposition \ref{propGamm} $(iii)$, we know that the $\cT(a)$-balls of radius $r_n$ are disjoint and that for all $\Upsilon\in\ccB_{a,r_n}$, there exists a unique $\cT(a)$-ball $\Gamma\in \ccB_{a,r_1}$ such that $\Upsilon\subset\Gamma$. Let us enumerate $\ccB_{a,r_1}$ as $\{\Gamma_i, 1\leq i\leq Z_{a,r_1}\}$, and set
\begin{align*}
\forall  i \! \in \! \{1\ldots Z_{a,r_1}\}, \; \ccB^{(i)}_{a,r_n}\! & = \! \left\{\Upsilon \! \in\!  \ccB_{a,r_n} : \Upsilon\! \subset \! \Gamma_i\right\} \\
\; {\rm and}\quad & S^{(i)}_{a,\br,\beps}\! =\#\left\{\Upsilon\in\ccB^{(i)}_{a,r_n} : \textrm{$\Upsilon$ is 
$(\br,\beps)$-small}\right\}.
\end{align*}
One has
\begin{equation}\label{decompor1}
S_{a,\br,\beps}-\mu(\br,\beps)\la(\cT)=\sum\limits_{i=1}^{Z_{a,r_1}} \left(S^{(i)}_{a,\br,\beps}-\mu(\br,\beps)\la(\Gamma_i)\right)=:\sum\limits_{i=1}^{Z_{a,r_1}} X_i.
\end{equation}
Let us denote $b=a-r_1/2$ and  recall from (\ref{defGb}) the definition of the sigma-field $\cG_{b}$. Adapting the proof of Lemma \ref{Tabrownian}, it is not difficult to see that  under $\N_{b}$, conditionally on  $\cG_{b}$, 
and  conditionally on $\{Z_{a,r_1}=k\}$, the r.v. $X_1,\ldots X_k$ are independent and have the same law as \mbox{$S_{r_1/2,\br,\beps}-\mu(\br,\beps)\ell^{r_1/2}(\cT)$} under $\N_{r_1/2}$. Recalling from (\ref{Laplaceella})  that $\N\left[\ell^{r_1/2}(\cT)\right]=(2/r_1)\N_{r_1/2}\left[\ell^{r_1/2}(\cT)\right]=1$, we see that 
$$\N_b[X_1\mid\cG_b] =\N_{r_1/2}\left[S_{r_1/2,\br,\beps}-\mu(\br,\beps)\ell^{r_1/2}(\cT)\right]=0,$$
which explains the definition (\ref{defmu}).
We thus apply Lemma \ref{moment4} to get from (\ref{decompor1}):

\begin{equation}\label{majodecompor1}
\N_{b}\left[\left(S_{a,\br,\beps}-\mu(\br,\beps)\la(\cT)\right)^4\mid \cG_{b}\right] \leq 3\N_{r_1/2}\left[X_1^4\right]\N_{b}\left[Z_{a,r_1}^2\mid \cG_{b}\right].
\end{equation}
The second assertion in  Lemma \ref{moment4} entails  $\N_{r_1/2}[X_1^4]\leq 2\N_{r_1/2}\left[S_{r_1/2,\br,\veps}^4\right]$. Moreover, we can use that  $S_{r_1/2,\br,\beps}$ is smaller than $Z_{r_1/2}(r_n)$, the total number of $\cT(r_1/2)$-balls of radius $r_n$ which has under $\N_{r_1/2}$ a geometric distribution with success probability $r_n/r_1<1/2$. Thus, 
\begin{equation}\label{majorecentre}
\N_{r_1/2}\left[X_1^4\right]
\leq 2\N_{r_1/2}\left[S_{r_1/2,\br,\veps}^4\right]\leq 2\N_{r_1/2}\left[Z_{r_1/2}(r_n))^4\right]\leq \frac{48}{1-r_n/r_1}\left(\frac{r_1}{r_n}\right)^4\leq 96\left(\frac{r_1}{r_n}\right)^4.
\end{equation}

In addition, according to the branching property, under $\N_{b}$, conditionally on $ \cG_{b}$,  $Z_{a,r_1}$ is a Poisson variable with mean $\N\left(h(\cT)>r_1/2\right)\ell^{b}(\cT)=(2/r_1)\ell^{b}(\cT)$. Thus, 
\begin{equation}\label{majoZar1}
\N_{b}\left[Z_{a,r_1}^2\right]
=(2/r_1)\N_{b}\left[\ell^{b}(\cT)\right]+ (2/r_1)^2\N_{b}\left[\ell^{b}(\cT)^2\right].
\end{equation}
Recalling that  $b=r_1/2$, we get  $(2/r_1)\N_{b}\left[\ell^{b}(\cT)\right]=1$ and $(2/r_1)^2\N_{b}\left[\ell^{b}(\cT)^2\right]=8a^2/r_1^2$. We assumed that $a/r_1>1$, thus $\N_{b}\left[Z_{a,r_1}^2\right]\leq 9a^2/r_1^2$. 
Combined with  (\ref{majodecompor1}) and  (\ref{majorecentre}) it entails
$$
\N_{b}\left[\left(S_{a,\br,\beps}-\mu(\br,\beps)\la(\cT)\right)^4\right]\leq c_0a^2\frac{r_1^2}{r_n^4},
$$
with $c_0$ a positive constant smaller than $(1/2)10^4$. This implies the desired result, using that $\N\left(h(\cT)>b\right)=1/b\leq 2/a$. 
\cqfd

\bigskip

We state now the main technical Lemma of the paper. Let us recall from (\ref{defmu}) the definition of $\mu(\br,\beps)$. The proof of the lemma makes use of the spinal decomposition described in Section  \ref{basicfactsBrownian}. In particular, a geometric argument allows to rely the problem to the variables introduced in (\ref{Lstardef}). 

%
%
%

%
%

\begin{lem}\label{munu}
Let $\br=\{r_1,\ldots,r_n\}$, where $r_1>\ldots>r_n>0$, and $\beps=\left\{\veps_1,\ldots,\veps_{n\!-\!1}\right\}$, where $\veps_1>\ldots>\veps_{n\!-\!1}>0$. The following inequality holds :
\begin{equation}
\mu(\br,\beps)\leq \frac{5}{r_n}\sqrt{\prod\limits_{k=1}^{n-1}\P
\left(\Lambda^*_{r_{k+1}, r_{k}}\leq\veps_k\right)}.
\end{equation}
\end{lem}

\noindent{\bf Proof.}
Let $\br=\{r_1,\ldots,r_n\}$ and $\{\veps_1,\ldots\veps_{n\!-\!1}\}$ as above. In that proof, we denote, for convenience, $b=r_1/2$; hence, a dependency with respect to $b$ is actually a dependency with respect to $\br$.   Let us consider  $\Gamma$  a $\cT(b)$-ball of radius $r_n$ and recall the notation (\ref{enlargedball}).  The ball $\Gamma$ is $(\br,\beps)$-small iff (\reff{condballsmall}) holds. But, for all $\sigma\in\Gamma, k\in\llbracket 1, n\!-\!1\rrbracket$, $$\Gamma[r_k]=\Gamma(\sigma,r_k)\supset\Gamma(\sigma,r_k)\setminus\Gamma(\sigma,r_{k+1}).$$
Thus, if $\Gamma$ is $(\br,\beps)$-small, then all the vertices in $\Gamma$ belong to the set  

\begin{equation}\label{defsetsmall}
\mathscr{S}(\br,\beps):=\left\{\sigma\in\cT(b) :  \forall k\in\{1\ldots n\!-\!1\}\quad \lb\left(\Gamma(\sigma,r_k)\setminus\Gamma(\sigma,r_{k+1})\right)\leq\veps_k\right\}.
\end{equation}
The last set is easy to handle  using the spinal decomposition. Indeed, according to (\ref{spinalrings}) and  the independence stated  in Lemma \ref{otarieimperiale}, one has 
\begin{equation}\label{defnu}
\nu(\br,\beps):=\N\left[\int\lb(\ddr \sigma)\un_{\{\sigma\in\mathscr{S}(\br,\beps)\}}\right]=\prod\limits_{k=1}^{n-1}\P
\left(\Lambda^*_{r_{k+1}, r_{k}}\leq\veps_k\right)
\end{equation}
To rely $\mu(\br,\beps)$ and $\nu(\br,\beps)$, one can write

\begin{equation}\label{trick}
\un_{\{\Gamma{\rm\  is \ }(\br,\beps)-{\rm small}\}}\leq
 \un_{\{\lb(\Gamma)\leq r_n\sqrt{\nu(\br,\beps)}\}}\ +\ \frac{\lb(\Gamma)}{r_n\sqrt{\nu(\br,\beps)}} \un_{\{\Gamma{\rm\  is \ }(\br,\beps)-{\rm small}\}}.
\end{equation}

\noi  Moreover, (\ref{defsetsmall}) entails that   $\lb(\Gamma)\un_{\{\Gamma{\rm\  is \ }(\br,\beps)-{\rm small}\}}\leq\int_\Gamma \lb(\ddr\sigma)\un_{\{\sigma\in\mathscr{S}(\br,\beps)\}}$. Recall now from Proposition \ref{propGamm} $(i)$   that the balls of the set $\ccB_{b,r_n}$ are pairwise disjoint. Summing in (\ref{trick}) over this set entails

\begin{equation}\label{tricksum}
S_{b,\br,\beps}\leq 
\sum\limits_{\Gamma\in\ccB_{b,r_n}}\un_{\left\{\lb(\Gamma)\leq r_n\sqrt{\nu(\br,\beps)}\right\}}
\ +\ 
\frac{\int\lb(\ddr \sigma)\un_{\{\sigma\in\mathscr{S}(\br,\veps)\}}}{r_n\sqrt{\nu(\br,\beps)}}.
\end{equation}

\noi Now, recalling Lemma \ref{Tabrownian}, we compute
\begin{align*}
\N_{b}\left[\sum\limits_{\Gamma\in\ccB_{b,r_n}}\un_{\left\{\ell^{b}(\Gamma)\leq r_n\sqrt{\nu(\br,\beps)}\right\}}\right]
&=\N_{b}\left[Z_{b,r_n}\right]\left(1-\exp\left(-(2/r_n)r_n\sqrt{\nu(\br,\beps)}\right)\right)\\
&\leq \frac{r_1}{r_n}2\sqrt{\nu(\br,\beps)},
\end{align*}
so the $\N$-measure of the first term in (\ref{tricksum}) is smaller than  $\frac{2b^{\!-1}r_1}{r_n}\sqrt{\nu(\br,\beps)}$. Recalling that $b=r_1/2$, we get that the latter equals  $\frac{4}{r_n}\sqrt{\nu(\br,\beps)}$. 
Moreover, by the mere definition (\ref{defnu}), the $\N$-measure of  the second term in (\ref{tricksum}) equals $\frac{1}{r_n}\sqrt{\nu(\br,\beps)}$, so the first inequality is checked.
\cqfd

\section{Proof of Theorem \ref{mainth}.}\label{preuvetheoreme}

\noi The proof of Theorem \ref{mainth} will combine the following two theorems. 

\begin{thm}\label{upper}
Let $\kappa\in (\frac{1}{2},\infty)$ and $m\in(0,\frac{1}{2})$. Then, there exists a Borel subset $\bV=\bV(\kappa,m)\subset\bC^0$ such that $\N\left(\bC^0\setminus \bV\right)=0$ and such that 
$${\rm on\ \bV},\quad {\rm for\ all\ Borel\ subset\ }\mathcal{A}\subset\cT,\quad  \forall a\in[m,m^{\!-1}],\quad \la(\mathcal{A})\leq \kappa\cH_g\left(\mathcal{A}\cap\cT(a)\right).$$
\end{thm}

\noi For all $a,\alpha\in(0,\infty)$, let us set 
\begin{equation}\label{defDelta}
\Delta_a^\alpha:=\left\{
\sigma\in\cT(a) : \limsup\limits_{r\to 0} \frac{\la\left(B(\sigma,r)\right)}{g(r)}< \alpha
\right\}.
\end{equation}

\begin{thm}\label{BP}
Let $\alpha\in(0,\frac{1}{2})$ and $m\in(0,1/2)$. Then, there exists a Borel subset $\bV^\prime=\bV^\prime(\alpha,m)\subset\bC^0$ such that $\N\left(\bC^0\setminus \bV^\prime\right)=0$ and such that
$${\rm on\ \bV^\prime},\quad  \forall a\in[m,m^{\!-1}],\quad \cH_g\left(\Delta_a^\alpha\right)=0.$$
\end{thm}

The proofs of Theorem \ref{upper} and  \ref{BP} share a common strategy, taken from Perkins \cite{Per88,Per89}. We need to control the mass, or the number of "bad" $\cT(a)$-balls where "bad" means too large or too small. And we want to do it uniformly for all levels $a$. This problem will be linked with a discrete one using a finite grid, and the measure or the number of bad $\cT(a)$-balls will be compared with a convenient multiple of $\la(\cT)$, the total mass at level $a$.

\subsection{Proof of Theorem \ref{upper}.}
\subsubsection{Large balls.}
Let us fix a level $a\in(0,\infty)$, and recall from section \ref{basicfactsBrownian} the definition of the sigma-field $\cG_a$, generated by the excursion below level $a$. We also recall  the definition of $\cT(a)$-balls (\ref{Taboule}). We fix  a threshold $y\in (0,\infty)$ and we consider the following set of "large" points on $\cT(a)$ :
\begin{equation}\label{defheavyballs}
\cL_{a,r,y}=\{\sigma\in \cT(a) : \ell^a(\Gamma(\sigma,r))> y\}.
\end{equation}
According to  Lemma \ref{Tabrownian}, the "total large mass" $\la\left(\cL_{a,r,y}\right)=\sum_{\Gamma\in\ccB_{a,r}}\la(\Gamma)\un_{\{\la(\Gamma)>y\}}$ is $\cG_a$-measurable.

\begin{lem}\label{GDlarge}
For all $a,l,y,r,\D \in(0,\infty)$, for all $c\in(1,\infty)$,   
$$\N\Big( \la\big( \cL_{a,r,y/c}\big)\leq l \ ;\ \sup_{b\in[a,a+\D]} \ell^b\big( \cL_{b,r,y}\big)> 4l\Big) 
\leq \frac{1}{a}\exp\left(-l/\D\right)\ +\ \frac{2}{r}\exp\left(-(1\!-\!c^{\!-\!1/2})^2 y/\D\right).$$
\end{lem}

\noindent{\bf Proof.}
For all $a\in(0,\infty)$, and all $c\in(1,\infty)$, we  define $A_0$, a Borel subset of $\bC^0$, as the event

\begin{equation}\label{defA0}
A_0=\left\{\la\left(\cL_{a,r,y/c}\right)\leq l \ ;\ \sup\limits_{b\in[a,a+\D]}\ell^b\left(\cL_{b,r,y}\right)> 4l\right\}.
\end{equation}
We recall from Proposition \ref{propGamm} that for   $r\in(0,\infty)$, $\ccB_{a,r}=\{\Gamma_i, 1\leq i\leq Z_{a,r}\}$ is the collection of  $\cT(a)$-balls of radius $r$ at level $a$. For $\Gamma$ a $\cT(a)$-ball and $b\in[a,\infty)$, we defined  $\Gamma^b=\{\sigma\in\cT(b) : \ \exists\sigma^\prime\in\Gamma,  \sigma^\prime\in\llbracket\rho,\sigma\rrbracket\}$ as the set of vertices at level $b$ having an ancestor in $\Gamma$ (see (\ref{enumpostboules}) for details).
Next we define $A_1$ a Borel subset of $\bC^0$ as the event
\begin{equation}\label{omega1}
A_1:=\left\{\exists i\in\{ 1,\ldots, Z_{a,r}\},\quad \la(\Gamma_i)\leq y/c {\rm \ \ and \ } \sup\limits_{b\in[a,a+\D]}\ell^b\left(\Gamma_i^b\right)>y\right\},
\end{equation}
\noi and set

\begin{equation}\label{defdescLb}
\cL^b_{a,r,y/c}:=\bigcup\limits_{\stackrel{i\in\{ 1,\ldots, Z_{a,r}\}}{\la(\Gamma_i)\geq y/c}}\Gamma_i^b\quad\subset\cT(b),
\end{equation} 
which is the set of all vertices at level $b$ having an  "large" ancestor at level $a$. We prove the following :
\begin{equation}\label{inclusioninterpo}
{\rm on\ } \bC^0\setminus A_1,\quad  \forall b\in[a,a+\D] \quad \cL_{b,r,y}\subset\cL^b_{a,r,y/c}.
\end{equation}
\noi 
\textit{Proof of}(\ref{inclusioninterpo}). 
Let $b\geq a\geq r/2>0$ and let $\sigma\in\cL_{b,r,y}$. Thus, the ball $\Gamma:=\Gamma(\sigma,r)\in\ccB_{b,r}$ is such that $\lb(\Gamma)\geq y$. Let $\sigma_a$ the unique ancestor of $\sigma$ at level $a$, namely $d(\rho,\sigma_a)=a$ and $\sigma_a\in\llbracket \rho,\sigma\rrbracket$. We set $\Upsilon:=\Gamma(\sigma_a,r)\in\ccB_{a,r}$ and  we first claim that $\Gamma\subset\Upsilon^b$. Indeed, let $\sigma^\prime\in\Gamma$ (so $d(\sigma,\sigma^\prime)<r$) and let $\sigma^\prime_a$ the unique ancestor of $\sigma^\prime$ at level $a$. Recalling that $\sigma\wedge\sigma^\prime$ stands for the most recent common ancestor of $\sigma$ and $\sigma^\prime$,  one get 
$d(\rho,\sigma\wedge\sigma^\prime)=\frac{1}{2}\left(2b-d(\sigma,\sigma^\prime)\right)$.
Then two cases may occur. First, if $d(\sigma,\sigma^\prime)\leq 2(b-a)$, then $d(\rho,\sigma\wedge\sigma^\prime)\geq a$, thus $\sigma_a=\sigma'_a$ and $\sigma'\in\Upsilon^b$. If  $d(\sigma,\sigma^\prime)\in (2(b-a),r)$. Then, one has $d(\rho,\sigma\wedge\sigma^\prime)<a$. We deduce from that inequality that  $\sigma_a\neq\sigma^\prime_a$ and that  $\sigma_a\wedge\sigma'_a=\sigma\wedge\sigma^\prime$. Hence, 
\begin{align*}
d(\sigma_a,\sigma^\prime_a)&=2a-2d\left(\rho,\sigma_a\wedge\sigma^\prime_a\right)\\
&= 2b-2d\left(\rho,\sigma\wedge\sigma^\prime\right) + 2a-2b\\
&=d(\sigma,\sigma^\prime)-2(b-a)<r.
\end{align*}
Thus $\sigma_a^\prime\in\Upsilon$ and $\sigma^\prime\in\Upsilon^b$, which ends the proof of the inclusion $\Gamma\subset\Upsilon^b$. We then get $\lb(\Upsilon^b)\geq \lb(\Gamma)> y$. On $\bC^0\setminus A_1$, one cannot have both $\la(\Upsilon)\leq y/c$ and $\lb(\Upsilon^b)>y$, which entails that here, $\la(\Upsilon)>y/c$.\\
To sum up,  on $\bC^0\setminus A_1$, a  vertex $\sigma$, taken in $\cL_{b,r,y}$, has an ancestor in a ball $\Upsilon$, such that $\la(\Upsilon)\geq y/c$. Thus, this ancestor belongs to $\cL_{a,r,y/c}$ and $\sigma\in \cL^b_{a,r,y/c}$.\\
\textit{End of the proof of} (\ref{inclusioninterpo}).\\

Let us finish the proof of the lemma. From (\ref{inclusioninterpo}), we see that 
\begin{equation}\label{majoA12}
\N(A_0)\ \leq\   \N\left(A_1\right)\ +\ \N\Big(A_0\cap (\bC^0\setminus A_1)\Big)\leq \  \N\left(A_1\right)\ +\ \N\left( A_2\right),
\end{equation}
where $A_2$ is defined by 
\begin{equation}\label{A2}
A_2:=\left\{\la\left(\cL_{a,r,y/c}\right)\leq l \ ;\ \sup\limits_{b\in[a,a+\D]}\ell^b\left(\cL^b_{a,r,y/c}\right)\geq 4l\right\}.
\end{equation}

We control $\N(A_1)$ and $\N(A_2)$ thanks to Lemma $\ref{GDFeller}$ (ii). Indeed,  Lemma \ref{branchingRK} states that  under $\N_a$, conditionally on $\cG_a$, 
the processes $\big(\ell^{a+a^\prime}(\Gamma_i^{a+a^\prime}), a^\prime\geq 0\big)$,$1\leq i\leq Z_{a,r}$ are independent Feller diffusions started at $\la(\Gamma_i), 1\!\leq\! i\!\leq\! Z_{a,r}$. Since $\lb\big(\cL^b_{a,r,y/c}\big)=\sum_{i=1}^{Z_{a,r}}\un_{\{\la(\Gamma_i)>y/c\}}\lb(\Gamma_i^b) $, it implies that  $\big(\ell^{a+a^\prime}(\cL^{a+a^\prime}_{a,r,y/c}), a^\prime\geq 0\big)$ is a Feller diffusion started at $\la\big(\cL_{a,r,y/c}\big)$. Thus, on the one hand, sub-additivity and Lemma \ref{GDFeller} (ii) entails
\begin{align}\label{majoA1}
\N\left(A_1\right)
&\leq \frac{1}{a}\N_a\left[\sum\limits_{i=1}^{Z_{a,r}}\un_{\{\la(\Gamma_i)\leq  y/c\}}\exp\left(-\D^{\!-1}\left(\sqrt{y}-\sqrt{\la(\Gamma_i)}\right)^2\right)\right]\nonumber\\
&\leq \frac{1}{a}\exp\left(-(1\!-\!c^{\!-\!1/2})^2 \D^{\!-1}y\right)\N_a[Z_{a,r}]=\frac{2}{r}\exp\left(-(1\!-\!c^{\!-\!1/2})^2\D^{\!-1} y\right).
\end{align}
On the other hand, 
\begin{align}\label{majoA2}
\N\left(A_2\right)
&\leq \frac{1}{a}\N_a\left[\un_{\{\la\left(\cL_{a,r,y/c}\right)\leq l\}} \exp\left(-\D^{\!-1}\left(2\sqrt{l}-\sqrt{\la\left(\cL_{a,r,y/c}\right)}\right)^2
\right)\right]\nonumber\\
&\leq \frac{1}{a}\exp\left(-l/\D\right).
\end{align}
Hence, the desired result follows from  (\ref{majoA12}), (\ref{majoA1}), and  (\ref{majoA2}). 
\cqfd

\bigskip

Recall that $g(r)=r\log\log(1/r)$. 
We fix $\kappa\in(\frac{1}{2},\infty)$, and we shall   apply the previous lemma with  $y=\kappa g(r)$.  The next lemma allows to control $\la\left(\cL_{a,r,\kappa g(r)}\right)$ uniformly for all levels $a$. Its proof involves a discrete grid :  for  $m<1/2$ and $r\in(0,\infty)$, we set 
\begin{equation}\label{defGr}
G(r,m):=\Big\{m+k\D_r, k\in\mathbb{N}^*\Big\}\cap  [m,m^{\!-1}],
\end{equation}
where $\D_r$ is the mesh of the grid, defined by
\begin{equation}\label{defDr}
\D_r=r^{3/2}.
\end{equation} 
Note that $G(r,m)$ contains less than $(m\D_r)^{-\!1}$ points. 
\begin{lem}\label{interpoLarge}
Let $m\in(0,1/2)$. Let   $\kappa\in(\frac{1}{2},\infty)$ and $\beta\in(1,\infty)$ such that $2\kappa-\beta>0$. There exists a constant $r_1\in(0,\infty)$ only depending on $\kappa,\beta,m$, such that 
\begin{equation}\label{ineginterpoLarge}
\forall r\in(0,r_1), \quad \N\left(\sup\limits_{b\in[m,m^{\!-1}]} \lb\left(\cL_{b,r,\kappa g(r)}\right)>4\log(1/r)^{-\beta}\right)\leq \log(1/r)^{-2}.
\end{equation}
\end{lem}

\noindent{\bf Proof.}
In what follows, we denote $T_0$ the left-hand-side of (\ref{ineginterpoLarge}).
Let us consider $c\in(1,\infty)$ such that $2\kappa/c-\beta>0$.  Recall that $G(r,m)$ stands for the grid defined by (\ref{defGr}). Then we have $T_0\leq T_1+T_2$, where :  
\begin{align*}
T_1&=\N\left(\sup\limits_{a\in G(r,m)}\la\left(\cL_{a,r,\kappa g(r)/c}\right)\leq \log(1/r)^{-\beta} \ ;\ \sup\limits_{b\in[m,m^{\!-1}]}\ell^b\left(\cL_{b,r,\kappa g(r)}\right)\geq 4\log(1/r)^{-\beta}\right),\\
T_2&=\N\left(\sup\limits_{a\in G(r,m)}\la\left(\cL_{a,r,\kappa g(r)/c}\right)> \log(1/r)^{-\beta}\right).
\end{align*}

\noi Using sub-additivity and Lemma \ref{GDlarge}, one get 
\begin{align*}
T_1&= \N\Big( \bigcup_{a\in G(r,m)}\Big \{ \la\big( \cL_{a,r,\kappa g(r)/c}\big) \leq \log(1/r)^{-\beta} \ ;\ \sup\limits_{b\in[a,a+\D_r]}\ell^b\big( \cL_{b,r,\kappa g(r)}\big) \geq 4\log(1/r)^{-\beta}\Big\} \Big) \\
&\leq (m\D_r)^{\!-1}\!\!\! \!\! \sup\limits\limits_{a\in G(r,m)} \!\!\! \N\Big( \la\big(\cL_{a,r,\kappa g(r)/c}\big) \leq \log(1/r)^{-\beta} \ ;\ \sup\limits_{b\in[a,a+\D_r]} \!\!\!\! \ell^b\big( \cL_{b,r,\kappa g(r)}\big) \geq 4\log(1/r)^{-\beta}\Big) \\
&\leq (m\D_r)^{\!-1}\bigg(m^{\!-1}\exp\left(-\D_r^{\!-1}\log(1/r)^{-\beta}\right)\ +\ \frac{2}{r}\exp\left(-(1\!-\!c^{\!-\!1/2})^2\kappa\D_r^{\!-1} g(r)\right)\bigg)
\end{align*}
One has $\D_r^{\!-1}\log(1/r)^{-\beta}\geq r^{-1}$ and $\D_r^{\!-1}g(r)\geq r^{-1/2}$ for all $r$ sufficiently small. Thus, for example, $T_1\leq \exp(-r^{-\!1/4})\leq (1/2)\log(1/r)^{-2}$ for all $r$ sufficiently small. \\

Let us bound $T_2$. To that end, we set
\begin{equation}\label{deflambda}
\lambda(r,\kappa,c):=(2/r)\E\left[\mathcal{E}\un_{\{\mathcal{E}> \kappa g(r)/c\}}\right],
\end{equation}
where $\mathcal{E}$ is a r.v.~on $(\Omega,\cF, \P)$  exponentially distributed with mean $r/2$.
For fixed $\kappa$ and $c$, elementary computations entail 
\begin{align}\label{equivlambda}
\lambda(r,\kappa,c)&=(2/r)\P\left(\mathcal{E}>\kappa g(r)/c\right)\E[\kappa g(r)/c + \mathcal{E}]\nonumber\\
&=(2/r)\exp\left(-2(\kappa/c) \log\log 1/r\right)\left((\kappa/c)r\log\log 1/r + r/2\right)\nonumber\\
&\underset{r\to 0}{\sim} (2\kappa/c)\log(1/r)^{-2\kappa/c}\log\log 1/r.
\end{align}

\noi We set $T_2 \leq T_3 + T_4$, where 
\begin{align*}
T_3&=\N\left(\sup\limits_{a\in G(r,m)}\left\lvert\la\left(\cL_{a,r,\kappa g(r)/c}\right)-\lambda(r,\kappa,c)\la(\cT)\right\rvert > \frac{1}{2}\log(1/r)^{-\beta}\right),\\
T_4&=\N\left(\sup\limits_{a\in G(r,m)}\lambda(r,\kappa,c)\la(\cT)> \frac{1}{2}\log(1/r)^{-\beta}\right).
\end{align*}
By  sub-additivity and  a Markov inequality involving a moment of order $4$, we get 
\begin{align}\label{majoT3}
T_3&\leq (m\D_r)^{\!-1}\sup\limits_{a\in G(r,m)}\N\left(\left\lvert\la\left(\cL_{a,r,\kappa g(r)/c}\right)-\lambda(r,\kappa,c)\la(\cT)\right\rvert>\frac{1}{2}\log(1/r)^{-\beta}\right)\nonumber\\
&\leq (m\D_r)^{\!-1}2^4 \log(1/r)^{4\beta}\!\!\!\!\!  \sup\limits\limits_{a\in G(r,m)} \!\!\!\!  a^{\!-1} \N_a\left[\left(\la\big( \cL_{a,r,\kappa g(r)/c}\big) -\lambda(r,\kappa,c)\la(\cT)\right)^4\right].
\end{align}
Recall notation $\ccB_{a,r}=\{\Gamma_i\ ,1\!\leq\! i\! \leq \!Z_{a,r}\}$ for the set of $\cT(a)$-balls with radius $r$. Then, consider the decomposition 
$$\la\left(\cL_{a,r,\kappa g(r)/c}\right)-\lambda(r,\kappa,c)\la(\cT)=\sum\limits_{i=1}^{Z_{a,r}}X_i,$$
where $X_i:=\la(\Gamma_i)\left(\un_{\{\la(\Gamma_i)\geq \kappa g(r)/c\}}-\lambda(r,\kappa,c)\right)$.
Using Lemma \ref{Tabrownian}, we see that under $\N_a$, conditionally on $Z_{a,r}$, the random variables $\la(\Gamma_1),\ldots \la(\Gamma_{Z_{a,r}})$ are independent and exponentially distributed with mean $r/2$. Thus, the definition (\ref{deflambda}) of $\lambda(r,\kappa,c)$  entails that under  $\N_a$, conditionally on $Z_{a,r}$, the r.v. $X_1,\ldots X_{Z_{a,r}}$ are i.i.d., with mean $0$ and a moment of order $4$.  Then, by Lemma~\ref{moment4},
\begin{equation}\label{majolarge4}
\N_a\left[\left(\sum\limits_{i=1}^{Z_{a,r}}X_i\right)^{\!\!\!\! 4}\right]
\leq 3\N_a(X_1^4)\N_a\left[Z_{a,r}^2\right].
\end{equation}
From (\ref{equivlambda}), we know that $\lambda(r,\kappa,c)\overset{r\to 0}{\longrightarrow}0$,  so for all sufficiently small $r$, $\lambda(r,\kappa,c)\leq 1/2$ and $\lvert X_1^4\rvert \leq\la(\Gamma_1)$, which implies  $\N_a[X_1]\leq\N_a[\la(\Gamma_1)^4]=\frac{3}{2}r^4$  for all  sufficiently small $r$. Moreover, $Z_{a,r}$ is under $\N_a$ a geometric r.v. with "success" probability $p=r/2a$ (see Lemma \ref{Tabrownian}), thus $\N_a\left[Z_{a,r}^2\right]=(2-p)/p^2\leq 8a^2/r^2$.  Combining (\ref{majoT3}) and (\ref{majolarge4}), we get, for all sufficiently small $r$,
\begin{equation}
T_3\leq 3.2^4.(m\D_r)^{-\!1}\log(1/r)^{4\beta}\sup\limits_{a\in G(r,m)}a^{-\!1}\frac{3r^4}{2}\frac{8a^2}{r^2}\leq  10^3 m^{\!-2} \log(1/r)^{4\beta} r^{1/2},
\end{equation}
recalling that $\D_r=r^{3/2}$. Observe now  that the right hand side is smaller than $(1/4)\log(1/r)^{-2}$ for all  sufficiently small $r$.\\

For the term $T_4$,  Lemma \ref{GDella} entails 
\begin{align}
T_4&\leq \N\left(\sup\limits_{b\in [m,m^{-\!1}]}\lb\left(\cT\right)> \frac{1}{2}\lambda(r,\kappa,c)^{-\!1}\log(1/r)^{-\beta}\right)
\nonumber\\
&\leq (2/m)\exp\left(-(m/4)\lambda(r,\kappa,c)^{-\! 1}\log(1/r)^{-\!\beta}\right).
\end{align}
By (\ref{equivlambda}), $$\lambda(r,\kappa,c)^{-\!1}\log(1/r)^{-\beta}\underset{r\to 0}{\sim}\frac{c}{2\kappa}\log(1/r)^{2\kappa\!/c-\beta}\log\log(1/r)^{-\!1}.$$
Recall that  $2\kappa/c>0$ and take $\veps\in(0,2\kappa/c\!-\!\beta)$. Thus, for all sufficiently small $r$,
$$T_4\leq (2/m)\exp\left(-\log(1/r)^{\veps}\right),$$
which is smaller than $(1/4)\log(1/r)^{-2}$  for all sufficiently small $r$. 
\cqfd

\subsubsection{Proof of Theorem \ref{upper}.}
Let $\kappa\in(1/2,\infty)$, and let $m\in(0,1/2)$.   Let $\beta\in(1,\infty)$ such that $2\kappa-\beta>0$. For all $a\in(0,\infty)$, $y\in(1,\infty)$ recall from (\ref{defDelta}) the definition   :  
\begin{equation}
\Delta_a^{y\kappa}=\left\{\sigma\in\cT(a) : \limsup\limits_{r\to 0}\frac{\la(B(\sigma,r))}{g(r)}< y\kappa\right\}.
\end{equation}

\noi For any $p\in\bbN$, set  $r_p:=y^{-p}$. By Lemma \ref{interpoLarge},   for all sufficiently large $p$, 
\begin{equation}
\N\left(\sup\limits_{a \in[m,m^{\!-1}]} \la\left(\cL_{a,r_p,\kappa g(r_p)}\right)>4\log(1/r_p)^{-\beta}\right)\leq  \log(1/r_p)^{-\!2}=\log(y)^{-\!2}p^{-2},
\end{equation}
whose sum over $p$ is finite. By Borel Cantelli lemma,
\begin{equation}\label{BCella}
\textrm{$\N$-a.e., for all sufficiently large $p$,}\quad \sup\limits_{a \in[m,m^{\!-1}]} \la\left(\cL_{a,r_p,\kappa g(r_p)}\right)\leq 4\log(1/r_p)^{-\beta}.
\end{equation}
Moreover, $\log (1/r_p)^{-\!\beta}=\log(y)^{\!-\beta}p^{\!-\beta}$, and recall that $\beta>1$. Thus, (\ref{BCella}) entails that there exists a Borel  subset $\bV_y\subset \bC^0$, such that $\N(\bC^0\setminus\bV_y)=0$, and on $\bV_y$: 
$$  \forall a\in[m,m^{\!-1}],\; \quad 
\sum\limits_{p=1}^\infty \la\big( \cL_{a,r_p,\kappa g(r_p)}\big)=\sum\limits_{p=1}^\infty \la\left(\{ \sigma : \la(B(\sigma,r_p))\! >\! \kappa g(r_p) \}\right)\! <\! \infty. $$

We can apply again the  Borel-Cantelli Lemma, to the finite measures $\la$ to get that, 
\begin{equation}\label{limsuprp}
{\rm on \ }\bV_y\quad \forall a\in[m,m^{\!-1}],\quad \la(\ddr \sigma)\textrm{-a.e.}\quad \exists p_0(a,\sigma),\quad\forall p\geq p_0(a,\sigma), \quad \frac{\la(B(\sigma,r_p))}{g(r_p)}\leq \kappa .
\end{equation}

If $u\in(r_{p+1},r_p]$, one has $\frac{\la(B(\sigma,u))}{g(u)}< \frac{\la(B(\sigma,r_p))}{g(r_{p+1})}\leq y\frac{\la(B(\sigma,r_p))}{g(r_p)}$. Combined with (\ref{limsuprp}), this entails that on    $\bV_y$, for all $a$ in $[m,m^{\!-1}]$, for $\la$-almost every $\sigma$ in $\cT(a)$,  $\limsup_{r\to 0}\la(B(\sigma,r))/g(r)<y\kappa$. This can be rewritten in 

\begin{equation}\label{limsuprp2}
{\rm on\ } \bV_y,\quad \forall a\in[m,m^{\!-1}],\quad \la\left(\cT(a)\setminus\Delta_a^{y\kappa}\right)=0.
\end{equation}

Now set $\bV=\bigcap\{\bV_y; y>1; y\in\bbQ\}$. Clearly,  $\N(\bC^0\setminus\bV)=0$  and by monotonicity, for all $\kappa^\prime\in(\kappa,\infty)$,  $\cT(a)\setminus\Delta_a^{\kappa^\prime}\subset\bigcup\limits_{y>1; y\in\bbQ}\left\{\cT(a)\setminus\Delta_a^{y\kappa}\right\}$. It follows easily  from  (\ref{limsuprp2}) that 

\begin{equation}\label{limsuprp3}
{\rm on \ } \bV, \quad\forall a\in[m,m^{\!-1}],\quad\forall \kappa^\prime\in(\kappa,\infty)\quad  \la\left(\cT(a)\setminus\Delta_a^{\kappa^\prime}\right)=0.
\end{equation}
Thus, using Lemma \ref{Comparison1}, we get :
\begin{align*}
{\rm on \ } \bV\quad \forall \mathcal{A} &{\rm \ Borel\  subset\  of\  }\cT\quad  \forall a\in[m,m^{\!-1}]\quad \forall \kappa^\prime\in(\kappa,\infty)\\
&\la\left(\mathcal{A}\right)=\la\left(\mathcal{A}\cap\Delta_a^{\kappa^\prime}\right)\leq  \kappa^\prime \cH_g\left(\mathcal{A}\cap\Delta_a^{\kappa^\prime}\right)\leq \kappa^\prime\cH_g\left(\mathcal{A}\cap\cT(a)\right).
\end{align*}
This ends the proof of Theorem \ref{upper} letting $\kappa^\prime\searrow\kappa$.

\subsection{Proof of Theorem \ref{BP}.}

\subsubsection{Small balls.}
For given level $a\in(0,\infty)$ and $r\in(0,\infty)$ we recall the notation  $\ccB_{a,r}$ for the set of $\cT(a)$-balls of radius $r$. We recall from (\ref{enlargedball}) that for $r\geq r^\prime>0$, a ball $\Gamma\in\ccB_{a,r^\prime}$ is contained in a unique ball in $\ccB_{a,r}$, denoted $\Gamma[r]$. Let  $\br=\{r_1>\ldots>r_n\}$ and $\beps=\{\veps_1>\ldots>\veps_{n\!-\!1}\}$.  Recall from \reff{condballsmall} that  $\Gamma$, a $\cT(a)$-ball of radius $r_n$ is $(\br,\beps)$-small iff
$$\forall k\in\llbracket 1,n\!-\!1\rrbracket \quad\la\left(\Gamma[r_k]\right)\leq \veps_k.$$
The total number of $(\br,\beps)$-small balls at level $a$ is denoted by  $S_{a,\br,\beps}$ (see \reff{defZs}). For $u\in(0,\infty)$, we write $u\br$ for the set $\{ur_1>\ldots>ur_n\}$. We recall from (\ref{postboule}) the following notation : if $\Gamma$ is a $\cT(a)$-ball, then, for all $b\geq a$, $\Gamma^b$ is the subset of all the vertices in  $\cT(b)$ having an ancestor in $\Gamma$. Namely, $\Gamma^b=\{\sigma\in\cT(b),  \exists \sigma^\prime\in\Gamma : \sigma^\prime\in\llbracket\rho,\sigma\rrbracket \}$.

\begin{lem}\label{GDsmall}
Let $a,\D\in(0,\infty)$, and $n\geq 2$. Let $\br=\{r_1>\ldots>r_n\}$ and $\beps=\{\veps_1>\ldots> \veps_{n\!-\!1}\}$ and $c\in(1,\infty)$. Let $\alpha\in (0,1/2)$ and $\talpha\in(\alpha,1/2)$. If $\D<\frac{c\!-\!1}{2c}r_n$, then

$$\N\left(\sup\limits_{b\in[a,a+\D]}S_{b,\br,\alpha\beps}
>S_{a,c^{\!-\!1}\br,\talpha\beps}\right)\leq \frac{2n}{r_n}\exp\left(\!-\!\left(\sqrt{\talpha}\!-\!\sqrt{\alpha}\right)^{\!2}\veps_{n\!-\!1}/\D\right).$$
\end{lem}

\noindent{\bf Proof.}
Let us denote $B_0=\big\{\sup_{b\in[a,a+\D]}S_{b,\br,\alpha\beps}
>S_{a,c^{\!-\!1}\br,\talpha\beps}\big\}$. Next, we  define the event $B_1$ by
\begin{equation}\label{defOmega3}
B_1=\big\{\exists k \! \in \! \{1\ldots , n\!-1\!\}, \;  \exists \Gamma \! \in\! \ccB_{a,r_k/c} \; : \; \la(\Gamma)\geq\talpha\veps_k \; \;  {\rm and}\;\inf_{b\in[a,a+\D]} \!\!\! \! \lb(\Gamma^b)\! <\! \alpha\veps_k\big\}.
\end{equation}
We will prove that $B_0\subset B_1$, that is to say
\begin{equation}\label{inclusionB}
{\rm on\ }\bC^0\setminus B_1,\quad \sup\limits_{b\in[a,a+\D]}S_{b,\br,\alpha\beps}
\leq S_{a,c^{\!-\!1}\br,\talpha\beps}.
\end{equation}
\noi
\textit{Proof of}(\ref{inclusionB}). 
We work deterministically on $\bC^0\setminus B_1$. The inequality (\ref{inclusionB}) follows from the following claim. 
\begin{center}
\textit{For every $b\in [a,a+\D]$, for every  $\Gamma$  a  $\cT(b)$-ball of radius $r_n$ which is $(\br,\alpha\beps)$-small,} \\
\textit{ there exists $\Upsilon$ a $\cT(a)$-ball of radius $r_n/c$ such that $\Upsilon$ is $(c^{\!-\!1}\br,\talpha\beps)$-small \textbf{and} $\Upsilon^b\subset\Gamma$.}
\end{center}
Assume that the latter is true. Then, to any  $(\br,\alpha\beps)$-small ball at level $b$ corresponds a $(c^{\!-\!1}\br,\talpha\beps)$-small ball at level $a$ and the correspondence is injective. Summing over all $\cT(b)$-ball, we obtain~(\ref{inclusionB}).\\
Now let $b\in[a,a+\D]$ and $\Gamma\in\ccB_{b,r_n}$ such that $\Gamma$ is $(\br,\alpha\beps)$-small. Let $\sigma\in\Gamma$ and let $\sigma_a$ its unique ancestor at level $a$. Namely $\sigma_a\in\cT(a)$ and $\sigma_a\in\llbracket\rho,\sigma\rrbracket$.  We denote $\Upsilon=\Gamma(\sigma_a,r_n/c)\in\ccB_{a,r_n/c}$ the $\cT(a)$-ball of radius $r_n/c$ that contains $\sigma_a$. We claim that $\Upsilon$ is $(c^{\!-\!1}\br,\talpha\beps)$-small and that $\Upsilon^b\subset\Gamma$. To prove this, we show
\begin{equation}\label{UpsinclusGam}
\forall k\in\{1,\ldots, n\}\ \quad \left(\Upsilon[r_k/c]\right)^b\subset\Gamma[r_k].
\end{equation}
Let $k\in\{1\ldots n\}$ and let $\gamma\in\left(\Upsilon[r_k/c]\right)^b$. Its unique ancestor at level $a$, denoted $\gamma_a$, is such that $\gamma_a\in\Upsilon[r_k/c]$.  Two cases may occur. First, if $d(\sigma,\gamma)\leq 2(b-a)$, then we have $2(b-a)\leq 2\D<\frac{c\!-\!1}{c}r_n<r_n\leq r_k$. The other case corresponds to $d(\sigma,\gamma)> 2(b-a)$. Then $d(\rho,\sigma\wedge\gamma)=\frac{1}{2}(2b-d(\sigma,\gamma))<a$. Thus, $\sigma\wedge\gamma=\sigma_a\wedge\gamma_a$ and we have 
\begin{align*}
d(\sigma,\gamma)&=2b-2d(\rho,\sigma\wedge\gamma)\\
&=2a-2d(\rho,\sigma_a\wedge\gamma_a)+2b-2a\\
&\leq d(\sigma_a,\gamma_a)+2\delta\\
&<\frac{r_k}{c} \ +\ \frac{c-1}{c}r_n\leq r_k,
\end{align*}
where we used in  the last line that  $\sigma_a\in\Gamma\subset \Gamma[r_k/c]$. In both cases, $d(\sigma,\gamma)<r_k$ so $\gamma\in\Gamma(\sigma,r_k)=\Gamma[r_k]$,
the last equality being a consequence of Proposition \ref{propGamm} $(ii)$, and the definition of $\Gamma=\Gamma(\sigma,r_n)$. Thus, (\ref{UpsinclusGam}) is proved and it implies 
$$\forall k\in\{1\ldots n\!-\!1\}\quad \lb\left(\left(\Upsilon[r_k/c]\right)^b\right)\leq\la\left(\Gamma[r_k]\right)\leq \alpha\veps_k,$$
which, on $\bC^0\setminus B_1$, implies 
$$\forall k\in\{1\ldots n\!-\!1\}\quad 
\la\left(\Upsilon[r_k/c]\right)\leq \talpha\veps_k.$$
This entails that $\Upsilon$ is $(c^{\!-\!1}\br,\talpha\beps)$-small. The inclusion $\Upsilon^b\subset\Gamma$ was proved at line (\ref{UpsinclusGam}) with $k=n$ because $\Upsilon=\Upsilon[r_n/c]\subset\Gamma[r_n]=\Gamma$.\\
\noi
\textit{End of the proof of} (\ref{inclusionB})\\

As in the proof of Lemma \ref{GDlarge}, we can use the fact that under $\N_a$, conditionally on  $\cG_a$, if  $\Gamma$ is a $\cT(a)$-ball, then the process $\left\{\ell^{a+a^\prime}(\Gamma^{a+a^\prime}), a^\prime\geq 0\right\}$ is a Feller diffusion started at $\la(\Gamma)$. Using sub-additivity and Lemma \ref{GDFeller} $(i)$, we get 
\begin{align}\label{majoB1}
\N\left(B_1\right)
&\leq\sum\limits_{k=1}^{n\!-\!1}\frac{1}{a}
\N_a\left[
\sum\limits_{i=1}^{Z_{a,r_k}}
\un_{\{\la(\Gamma_i)\geq\talpha\veps_k\}}
\exp\left(-\D^{\!-\!1}\left(\sqrt{\la(\Gamma_i)}-\sqrt{\alpha\veps_k}\right)^2\right)\right]\\
&\leq\frac{1}{a}\exp\left(-\D^{\!-\!1}\left(\sqrt{\talpha\veps_k}-\sqrt{\alpha\veps_k}\right)^2\right)\sum\limits_{k=1}^{n\!-\!1}\E\left[Z_{a,r_k}\right]
\end{align}
The proof is completed recalling that for all $k\in\{1\ldots n\!-\!1\}$, $\veps_k\leq \veps_{n\!-\!1}$, and that, by Lemma~\ref{Tabrownian}, \mbox{$\N_a\left[Z_{a,r_k}\right]=\frac{2a}{r_k}\leq\frac{2a}{r_n}$}. 
\cqfd

\vskip 0.2cm
Let us introduce
\begin{equation}\label{lesradis}
\forall j\in\mathbb{N},\quad r_j=2^{-j}\quad \textrm{and}\quad\veps_j=g(r_j)
\end{equation}
and then 
\begin{equation}\label{lesbonsentiers}
\forall p\in\mathbb{N},\quad j_p\! =\! \lfloor(4/3)^p\rfloor,\quad  \Rp \! =\! \{r_j, j_p\leq j\leq j_{p+1}\!-\!1\}\quad\textrm{and}\quad  \bepsp=\{\veps_j \, ; \,  j_p \! \leq \!  j< \! j_{p+1}\!-\!1\}.
\end{equation}

Let  $m\in(0,1/2)$, we also introduce the following discrete grid   
\begin{equation}\label{defG'p}
G^\prime(p,m):=\left\{m^{\!-1}+k\D_p, k\in\mathbb{N}^*\right\}\cap[m,m^{\!-1}],
\end{equation}
where $\D_p$ is the mesh of the grid, given by
\begin{equation}\label{defDp}
\D_p=r(j_{p+1})^{5/4}.
\end{equation} 
Note that $G^\prime(p,m)$ contains less than $(m\D_p)^{\!-1}$ points.

\begin{lem}\label{interposmall}
Let $\alpha\in(0,1/2)$, $m\in(0,1/2)$. For $p\in\mathbb{N}$, denote $u_p:=g\left(r(j_{p+1})\right)^{-1}p^{-2}$. Then there exists $p_0\in\mathbb{N}$ only depending on $\alpha,m$ such that for all $p\geq p_0$, 
\begin{equation}\label{ineginterposmall}
\N\left(\sup\limits_{b\in[m,m^{-\!1}]}S_{b,\Rp,\alpha\bepsp}>u_p\right)\leq p^{-2}.
\end{equation}
\end{lem}

\noindent{\bf Proof.}
Let $\talpha\in(\alpha,1/2)$ and  $c$ in $(1,\infty)$ such that $2c\talpha\in(0,1)$. In what follows, we denote $T'_0$ the left-hand-side of (\ref{ineginterposmall}). Observe that $T'_0\leq T'_1+T'_2$, where we have set
\begin{align*}
T'_1&=\N\left(\sup\limits_{a\in G'(p,m)}S_{a,c^{\!-\!1}\Rp\!,\talpha\bepsp}\leq u_p\ ; 
\sup\limits_{b\in[m,m^{-\!1}]}\!S_{b,\Rp,\alpha\bepsp}>u_p\right),\\
T'_2&=\N\left(\sup\limits_{a\in G'(p,m)}S_{a,c^{\!-\!1}\Rp\!,\talpha\bepsp}> u_p\right).
\end{align*}
Using sub-additivity and Lemma \ref{GDsmall}, we get 
\begin{align*}
T'_1&\leq \N\left(\bigcup\limits_{a\in G'(p,m)}\!\!\left\{\sup\limits_{b\in[a,a+\D_p]}\!S_{b,\Rp,\alpha\bepsp}>S_{a,c^{\!-\!1}\Rp\!,\talpha\bepsp}\right\}\right)\\
&\leq (m\D_p)^{\!-\!1}\sup\limits_{a\in G'(p,m)}\N\left(\sup\limits_{b\in[a,a+\D_p]}\!S_{b,\Rp,\alpha\bepsp}>S_{a,c^{\!-\!1}\Rp\!,\talpha\bepsp}\right)\\
&\leq (m\D_p)^{\!-\!1}\frac{2(j_{p+1}-j_p)}{r(j_{p+1})}\exp\left(\!-\!\left(\sqrt{\talpha}\!-\!\sqrt{\alpha}\right)^{\!2}\D_p^{\!-\!1}g(r(j_{p+1}\!-\!2))\right).
\end{align*}
One has $\D_p^{\!-\!1}g(r(j_{p+1}\!-\!2))\geq\D_p^{\!-\!1}g(r(j_{p+1})) =r(j_{p+1})^{\!-\!1/4}\log\log 1/r(j_{p+1})$, which implies that $T'_1$ is smaller than $(1/2)p^{-2}$, for all $p$ sufficiently large (it is obviously not a sharp bound). \\

Recalling the definitions \reff{defmu}, we set 

\begin{equation}\label{defmup}
\mu_p=\mu(c^{\!-1}\Rp\!,\talpha\bepsp)=\N\left(S_{r(j_p)/(2c),c^{\!-\!1}\Rp,\talpha\bepsp}\right).
\end{equation}
We will prove that $T'_2\leq T'_3+T'_4$, where 
\begin{align*}
T'_3&=\N\left(\sup\limits_{a\in G'(p,m)}\lvert S_{a,c^{\!-\!1}\Rp,\talpha\bepsp}-\mu_p\la\left(\cT\right)\rvert> u_p/2\right),\\
T'_4&=\N\left(\sup\limits_{a\in G'(p,m)}\mu_p\la\left(\cT\right)> u_p/2\right).
\end{align*}
By sub-additivity and a Markov inequality involving a moment of order 4, we get 
\begin{align}\label{borneT'3}
T'_3&\leq (m\D_p)^{\!-\!1}\sup\limits_{a\in G'(p,m)}\N\left(\lvert S_{a,c^{\!-\!1}\Rp,\talpha\bepsp}-\mu_p\la\left(\cT\right)\rvert> u_p/2\right)\nonumber\\
&\leq (m\D_p)^{\!-\!1} 2^4 u_p^{\!-\!4}\sup\limits_{a\in G'(p,m)}\N\left[\left(S_{a,c^{\!-\!1}\Rp,\talpha\bepsp}-\mu_p\la\left(\cT\right)\right)^{\!4}\right].
\end{align}
We want to apply Lemma \ref{Z4} with $\br=,c^{\!-\!1}\Rp$ and $\beps=\talpha\bepsp$. Thus, recalling (\ref{lesradis}) and (\ref{lesbonsentiers}), we check that for all sufficiently large $p$, $m/r(j_p)>1$ and $r(j_p)/r(j_{p+1}\!-\!1)$. Recalling that $c_0\in(0,10^4]$ is the universal constant given by Lemma \ref{Z4}, we get from (\ref{borneT'3})
\begin{equation}\label{borne2T'3}
T'_3\leq  (m\D_p)^{\!-\!1} 2^4 u_p^{\!-\!4}\sup\limits_{a\in G'(p,m)}
c_0a\frac{r(j_p)^2}{r(j_{p+1}\!-\!1)^4}\leq 2^4c_0 m^{\!-\!2} \frac{r(j_p)^2}{\D_p u_p^{4}r(j_{p+1})^4}.
\end{equation}
Recall that $u_p=g(r(j_{p+1}))^{\!-1}p^{\!-2}$, and by (\ref{lesradis}) and (\ref{lesbonsentiers}), we get $\log\log(1/r(j_p))\underset{p\to \infty}{\sim} p\log(4/3)$. Hence, $u_p\geq p^{-3}r(j_{p+1})^{-1}$ and (\ref{borne2T'3}) implies
\begin{equation}
T'_3\leq 2^4c_0 m^{\!-2}p^{12}\frac{r(j_p)^2\cancel{r(j_{p+1})^{4}}}{r(j_{p+1})^{5/4}\cancel{r(j_{p+1})^{4}}}.
\end{equation}
Now, one can plainly check that $\frac{r(j_p)^2}{r(j_{p+1})^{5/4}}$ is smaller than $r(j_p)^{1/3}$. Thus, $T'_3$ is smaller than $(1/4)p^{-\!2}$ for all $p$ sufficiently large.\\

For the  term $T'_4$, we use Lemma \ref{GDella} to obtain 
\begin{equation}\label{majoT'3}
T'_4\leq (2/m)\exp\left(-(m/4)u_p\mu_p^{-1}\right).
\end{equation}
Recalling (\ref{defmup}) and Lemma  \ref{munu}, we get that for all $p$, 
\begin{equation}
\mu_p\leq \frac{5}{r(j_{p+1})}\left(\prod\limits_{j=j_p}^{j_{p+1}-2}
\P\bigg(\Lambda^*_{r_{j+1}/c,r_j/c}\leq \talpha r_j\log\log(1/r_j)\bigg)\right)^{1/2}
\end{equation}
We want to get an lower bound of $u_p\mu_p^{\!-1}$, so we compute an upper bound for  $u_p^{\!-1}\mu_p$. Recalling that $u_p\geq p^{-3}r(j_{p+1})^{-\!1}$, one has 
\begin{equation}\label{upmup}
u_p^{\!-1}\mu_p\leq 5p^3\exp\left(\ \frac{1}{2}\sum\limits_{j=j_p}^{j_{p+1}-2}\log\left(1-q_j\right)\right)\leq 5p^3\exp\left(-\frac{1}{2} \sum\limits_{j=j_p}^{j_{p+1}-2}q_j\right),
\end{equation}

%
\noi where $q_j=\P\big(\Lambda^*_{r_{j+1}/c,r_j/c}> \talpha r_j\log\log(1/r_j)\big)$. Recalling that $r_j=2^{-j}$, it follows from Lemma \ref{otarieimperiale} that
\begin{align*}\label{evalqj}
q_j&=\left(1\!-\!\frac{1}{2}\right)^2\frac{2\talpha\cancel{r_j}\log\log 1/r_j}{\cancel{r_j}/c}\exp\left(-\frac{2\talpha\cancel{r_j}\log\log 1/r_j}{\cancel{r_j}/c}\right)+\left(1\!-\!\frac{1}{4}\right)\exp\left(-\frac{2\talpha\cancel{r_j}\log\log 1/r_j}{\cancel{r_j}/c}\right)\\
&\underset{j\to\infty}{\sim}\frac{\talpha c}{2}\log\log(1/r_j)e^{-2\talpha c\log\log(1/r_j)}\\
&\underset{j\to\infty}{\sim} c' \log(j)j^{-2\talpha c},
\end{align*}
\noi where $c'$ is a positive constant depending on $\alpha,\talpha,c$. We stress that the particular choice  of $c$ was made to ensure that $\chi:=1-2\talpha c$ is strictly positive, so that  the following is true for all large $p$ :
$$\sum_{j_p}^{j_{p+1}\!-2\!}q_j\geq \sum_{j_p}^{j_{p+1}\!-2\!}j^{-2\talpha c}\geq \int_{j_p}^{j_{p+1}\!-\!1}x^{-2\talpha c}\ddr x\underset{p\to\infty}{\sim} \chi^{-1}\left((4/3)^{\chi}\!-\!1\right)\left(\frac{4}{3}\right)^{p\chi}.$$

Thus, for all $p$ sufficiently  large, $\sum_{j_p}^{j_{p+1}}q_j\geq 2p$  which, combined  with (\ref{upmup}), entails that 
$u_p^{\!-1}\mu_p\leq 5p^3\exp\left(-p\right)$. Thus, $u_p\mu_p^{\!-1}\geq 5^{\!-1}p^{\!-3}e^p$. Finally, we see from (\ref{majoT'3}) that  $T'_3$ is smaller than $(1/4)p^{-2}$ for all $p$ sufficiently large, which ends the proof. \cqfd

\subsubsection{Proof of Theorem \ref{BP}.}
Let $\alpha\in(0,1/2)$. For a level $a\in(0,\infty)$,  we recall the definition (\ref{defDelta}) of $\Delta_a^\alpha$. To show that the $g$-Hausdorff measure  of $\Delta_a^\alpha$ is null, we need an efficient covering of this set. 
Let us recall the integer sequence $j_p=\lfloor(4/3)^p\rfloor$ and the radii $r_j=2^{-j}$. For $p\in\mathbb{N}$, we recall the definition of the finite subsets $\Rp=\{r_j, j_p\leq j\leq j_{p+1}\!-\!1\}$, and  $\bepsp=\{\veps_j, j_p\leq j<j_{p+1}\!-\!1\}$ where $\veps_j=g(r_j)$. 
Recalling the definition \reff{condballsmall} for small balls,  we set $$\mathscr{C}_n:=\bigcup\limits_{p= n}^{\infty}
\left\{\Gamma\in\ccB_{a,r(j_{p+1})} :\quad  \Gamma {\rm\  is\ }(\Rp,\alpha\bepsp){\rm -small}\right\}.$$
Observe that if $\sigma\in\Delta_a^\alpha$, then the $\cT(a)$-ball $\Gamma\left(\sigma,r(j_{p+1})\right)$ is  $(\Rp,\alpha\bepsp)$-small for all large $p$, thus for all $n\in\mathbb{N}$, we have $\Delta_a^\alpha\subset\mathscr{C}_n$. Let us recall the definition  (\ref{defHausdorff}) of Hausdorff measures, and the fact that the diameter of a $\cT(a)$-ball is smaller than its radius. We get 
\begin{equation}\label{majoLambda}
\forall a\in[m,m^{-\!1}]\quad\cH_g^{_{(r(j_{p+1}))}}\!\left(\Delta_a^\alpha\right)\leq \sum
\limits_{p=n}^{\infty}S_{a,\Rp,\alpha\bepsp}.\
g\left(r(j_{p+1})\right),
\end{equation}
because $\Delta_a^\alpha\subset\mathscr{C}_n$. Thus,
\begin{equation}\label{majoLambda2}
\forall a\in[m,m^{-\!1}]\quad\cH_g\left(\Delta_a^\alpha\right)\leq \limsup\limits_{n\to\infty}\sum
\limits_{p=n}^{\infty}S_{a,\Rp,\alpha\bepsp}.
g\left(r(j_{p+1})\right).
\end{equation}

\noi Now, let $m\in(0,1/2)$. Applying Lemma \ref{interposmall}, we easily get that 
$$\sum\limits_{p=1}^\infty\N\left(\sup\limits_{a\in[m,m^{-\!1}]}S_{a,\Rp,\alpha\bepsp}>u_p\right)<\infty,$$
where we recall the notation $u_p=g\left(r(j_{p+1})\right)^{-1}p^{-2}$. By Borel-Cantelli lemma  there exists a subset $\bV^\prime\subset\bC^0$ such that $\N\left(\bC^0\setminus\bV^\prime\right)=0$ and such that 
$${\rm on\ } \bV^\prime,\quad  g\left(r(j_{p+1})\right)\sup\limits_{a\in[m,m^{\!-1}]}S_{a,\Rp,\alpha\bepsp}\leq p^{-2},\quad {\rm for\ all\ suff.\ large\ }p.$$

\noi Combined with \reff{majoLambda}, we deduce on $\bV^\prime$, for  $a\in[m,m^{-\!1}]$, one has 
$$\cH_g\left(\Delta_a^\alpha\right)\leq \lim_{n\to\infty}\sum_{p=n}^{\infty}p^{-2}=0,$$ 
which is the desired result.

\subsection{Proof of Theorem \ref{mainth}.}
Let $\kappa\in(\frac{1}{2},\infty)$, $\alpha\in(0,\frac{1}{2})$, and $m\in(0,1/2)$. Theorem \ref{upper} entails that there exists a Borel subset $\bV=\bV(\kappa,m)\subset\bC^0$ such that $\N\left(\bC^0\setminus \bV\right)=0$ and 
\begin{equation}\label{final1}
{\rm on\ }\bV(\kappa,m),\quad {\rm for\ all\ Borel\ subset\ }\mathcal{A}\subset\cT,\quad  \forall a\in[m,m^{\!-1}],\quad \la(\mathcal{A})\leq \kappa\cH_g\left(\mathcal{A}\cap\cT(a)\right).
\end{equation}
Now, let us rewrite the definition (\ref{defDelta})
\begin{equation}
\Delta_a^\alpha=\left\{
\sigma\in\cT(a) : \limsup\limits_{r\to 0} \frac{\la\left(B(\sigma,r)\right)}{g(r)}< \alpha
\right\}.
\end{equation}
According to Theorem \ref{BP}, there exists a Borel subset $\bV^\prime=\bV^\prime(\alpha,m)\subset\bC^0$ such that $\N\left(\bC^0\setminus \bV^\prime\right)=0$ and
\begin{equation}\label{final2}
{\rm on\ \bV^\prime}(\alpha,m)\quad\forall a\in[m,m^{-\!1}]\quad \cH_g\left(\Delta_a^\alpha\right)=0.
\end{equation}
Let $\alpha'<\alpha$ and notice that $\cT(a)\setminus \Delta_a^\alpha\subset\left\{\sigma : \limsup\limits_{r\to 0}\frac{\la\left(B(\sigma,r)\right)}{g(r)}> \alpha'
\right\}$. Moreover, from (\ref{supportla}), we know that $\N$-a.e. for all $a\in(0,\infty)$, $\la(\cT\setminus\cT(a))=0$. 
Thus, on $\bV^\prime$, for all Borel subset $\mathcal{A}\subset\cT$, and for all $a\in[m,m^{-\!1}]$ and all $\talpha<\alpha$, Lemma \ref{Comparison2}  entails 
\begin{equation}\label{final3}
\la\left(\mathcal{A}\right)\geq \la\left(\mathcal{A}\cap(\cT(a)\setminus\Delta_a^\alpha)\right)
\geq \alpha'\cH_g\left(\mathcal{A}\cap(\cT(a)\setminus\Delta_a^\alpha)\right)=\alpha'\cH_g\left(\mathcal{A}\cap\cT(a)\right),
\end{equation}
where we used  (\ref{final2}) for the last equality. Letting $\alpha'\to\alpha$, we get 
\begin{equation}\label{final4}
{\rm on\ \bV^\prime}(\alpha,m)\quad{\rm for\ all\ Borel\ subset\ }\mathcal{A}\subset\cT,\quad\forall a\in[m,m^{-\!1}]\quad\la\left(\mathcal{A}\right)\geq\alpha\cH^g\left(\mathcal{A}\cap\cT(a)\right).
\end{equation}

Now, let us set
\begin{equation}\label{finalevent}
\tilde{\bV}:=\left(\bigcap\limits_{\stackrel{\kappa\in(1/2,\infty)\cap\bbQ}{_{m\in(0,1/2)\cap\bbQ}}}\!\!\bV(\kappa,m)\right)
\bigcap
\left(\bigcap\limits_{\stackrel{\alpha\in(0,1/2)\cap\bbQ}{_{m\in(0,1/2)\cap\bbQ}}}\!\!\bV^\prime(\alpha,m)\right).
\end{equation}
Clearly, $\tilde{\bV}$ is a Borel subset of $\bC^0$ such that  $\N\left(\bC^0\setminus \tilde{\bV}\right)=0$. Moreover,  combining \reff{final1} and \reff{final4}, we get that  on $\tilde{\bV}$, for all Borel subset $\mathcal{A}\subset\cT$, and for all level $a\in(0,\infty)$, one has $\la\left(\mathcal{A}\right)=\frac{1}{2}\cH_g\left(\mathcal{A}\cap\cT(a)\right)$.

\end{document}